\documentclass[11pt]{article}
\usepackage{amssymb}
\usepackage{amsmath}
\usepackage[dvips]{epsfig}
\usepackage{graphicx}
\usepackage{graphicx}
\usepackage[colorlinks=true,linkcolor=blue,citecolor=red]{hyperref}

\newtheorem{Lemma}{Lemma}[section]
\newtheorem{Theorem}{Theorem}

\newtheorem{Corollary}[Lemma]{Corollary}

\newtheorem{Definition}[Lemma]{Definition}

\newenvironment{Proof}
 {\begin{trivlist} \item[]{\bf Proof. }}
 {\QED\end{trivlist}}

\newenvironment{Acknowledgment}
 {\begin{trivlist}\item[]\textbf{Acknowledgments }}{\end{trivlist}}

\newcommand{\mcaption}[1]{\vspace*{-1ex}\begin{quote}\caption{\small #1}\end{quote}\vspace*{-5ex}}

\makeatletter\@addtoreset{figure}{section}\makeatother

\makeatletter\@addtoreset{equation}{section}\makeatother

\newcommand{\C}{\mathbb{C}}
\newcommand{\N}{\mathbb{N}}
\newcommand{\R}{\mathbb{R}}
\newcommand{\Z}{\mathbb{Z}}

\def\Re{\mathop{\mathrm{Re}}}

\newcommand{\per}{\mathrm{per}}
\newcommand{\loc}{\mathrm{loc}}

\renewcommand{\phi}{\varphi}

\def\QED{\mbox{}\hfill$\Box$}

\renewcommand{\leq}{\leqslant}
\renewcommand{\geq}{\geqslant}

\begin{document}

\thispagestyle{empty}

\vspace*{4ex}

\begin{center}
{\bf\LARGE Transverse spectral stability of small periodic\\[0.5ex]
traveling waves for the KP equation}
\\[4ex]
{\large\bf Mariana  Haragus}
\\
Laboratoire de Math\'ematiques,
Universit\'e de Franche-Comt\'e\\
16 route de Gray, 25030 Besan\c con, France\\
mharagus@univ-fcomte.fr
\end{center}

\vspace*{4ex}

\begin{abstract}
The Kadomtsev-Petviashvili (KP) equation possesses a four-parameter
family of one-dimensional periodic traveling waves. We study the
spectral stability of the waves with small amplitude with respect to
two-dimensional perturbations which are either periodic in the
direction of propagation, with the same period as the 
one-dimensional traveling wave, or non-periodic (localized or bounded). 
We focus
on the so-called KP-I equation (positive dispersion case), for which
we show that these periodic waves are unstable with respect to both types
of perturbations. Finally, we briefly discuss the KP-II equation, for
which we show that these periodic waves are spectrally stable with
respect to perturbations which are periodic in the
direction of propagation, and have long wavelengths in the transverse
direction. 
\end{abstract}

\vspace*{4ex}

\noindent {\bf Running head:} {Periodic waves of the KP
equation}


\vspace*{4ex}

\tableofcontents

\vspace*{4ex}

\begin{Acknowledgment}
This work was partially supported by the Agence Nationale de la Recherche
through grant ANR PREFERED.
\end{Acknowledgment}

\vspace*{4ex}

\section{Introduction}
\label{s:1}

We consider the Kadomtsev-Petviashvili (KP) equation 
\begin{equation}\label{e:kp}
\left( u_t-u_{xxx}-uu_x\right)_x+\sigma u_{yy}=0,
\end{equation}
in which $u(x,y,t)$ depends upon the spatial variables $x,y\in\R$, and
the temporal variable $t\in\R$, and $\sigma$ is equal to either $1$ or $-1$. 
This equation is a generalization to two spatial dimensions of the
well-known Korteweg-de Vries (KdV) equation
\begin{equation}\label{e:kdv}
u_t=u_{xxx}+uu_x,
\end{equation}
and has been derived as a model in the study of the transverse
stability of the line solitary waves, i.e., the stability of the
solitary waves of 
the KdV equation with respect to two-dimensional perturbations \cite{KP}.
In particular, it arises as a model for long waves 
in the water-wave problem, 
in the absence of surface tension when $\sigma=-1$, 
and in the presence of a strong surface tension
when $\sigma=1$. The equation with $\sigma=-1$
(negative dispersion)
is also called the KP-II equation, whereas the one with $\sigma=1$
(positive dispersion) is
called the KP-I equation. It turns out that the stability
properties of the line solitary waves strongly depend upon the sign of
$\sigma$: when $\sigma=-1$ the line solitary waves are transversely stable,
and when $\sigma=-1$ they are unstable \cite{APS,KP}.

In addition to the well-known solitary waves, the KdV equation possesses a
family of periodic traveling waves. In this paper we are interested in
the stability of these periodic waves as solutions of the KP equation,
i.e., their transverse stability. 
As solutions of the KdV equation, the stability of these 
periodic traveling waves is quite well-understood, though the results are
less complete than those for solitary waves. In contrast to the case
of localized solutions, for periodic waves of dispersive equations 
there are mainly to types
of stability results: nonlinear (orbital) stability with respect to
periodic perturbations which have the same period as the traveling
wave, and spectral stability with respect to non-periodic (localized
or bounded) perturbations.  The question of nonlinear stability with
respect to non-periodic perturbations is widely open, so far. 
For the KdV, it has been recently shown that the
periodic traveling waves are orbitally stable  with respect to
periodic perturbations \cite{ABS}, and spectrally 
stable with respect
to non-periodic perturbations \cite{BD,HK}.

Very recently, Johnson and Zumbrun \cite{JZ} derived an instability
criterion for the transverse instability of periodic waves for the
generalized KP equation, with respect to perturbations which are
periodic in the direction of propagation and have long wavelengths in
the transverse direction. In contrast, we restrict here to the case of
small periodic waves, and consider the question of spectral stability
for  more general perturbations for 
the KP-I equation, and for the same type of perturbations for the KP-II
equation. 
A particularity of the KP equation is that the spectral problem cannot
be directly formulated in terms of the spectrum of a linear operator,
because of the mixed time-space derivative $u_{tx}$ appearing in \eqref{e:kp}.
We present a functional set-up which leads to such a formulation, for
both periodic and non-periodic perturbations.
For periodic perturbations, 
we use the invertibility of the operator $\partial_x$ when
restricted to a space of periodic functions with zero mean, and for
non-periodic perturbations we rely  upon Floquet theory for
differential operators with periodic coefficients and the
invertibility of the operator $\partial_x+i\gamma$. For the analysis
of the spectra of the resulting operators we rely upon several tools from the 
perturbation theory for linear operators, which turn out to be
suitable for the analysis of spectra of differential operators with
periodic coefficients \cite{bbm,HLS}. In addition, we use 
the decomposition of these operators in the product of a skew-adjoint
with a self-adjoint operator, and a counting result which relates the
number of unstable eigenvalues of the product operator to the number
of negative eigenvalues of the self-adjoint operator \cite{HK}.
We focus on the KP-I equation, for  which
we show that the small periodic waves are transversely 
unstable, with respect to both types
of perturbations. 
For the KP-II equation, a complete analysis of the
spectra turns out to be more difficult, because of the different
behavior of the dispersion relation. Restricting to
perturbations which are periodic in the
direction of propagation and have long wavelengths in the transverse
direction, we show that the small periodic waves are spectrally
stable with respect to such perturbations. In particular, these
results show that the periodic waves have the same stability
properties as the solitary waves, but under the restriction to  long
wavelength transverse perturbations for the KP-II equation.

The paper is organized as follows. In Section~\ref{s:2}, 
we describe the family of the
one-dimensional periodic
traveling waves of the KP equation and give a parameterization of the
small waves. The spectral stability problem is
first formulated in Section~\ref{s:3}, and then reformulated for
periodic perturbations in Section~\ref{ss:specper}, and for
non-periodic perturbations in Section~\ref{ss:bloch}. 
The results in these latter two sections are presented for the
KP-I equation, but it is easily seen that the same results hold for the KP-II
equation. We discuss the spectra of the resulting linear
operators, and show the transverse spectral instability of the small
periodic waves of the KP-I equation for periodic perturbations in
Section~\ref{s:4}, and for non-periodic perturbations in
Section~\ref{s:5}. In Section~\ref{s:6}, we prove the 
transverse spectral stability result for the periodic waves of the
KP-II equation mentioned above, and we conclude with a brief
discussion in Section~\ref{s:7}.

\section{One-dimensional periodic traveling waves}
\label{s:2}

One-dimensional traveling waves of the KP equation \eqref{e:kp}
are solutions of the form 
\[u(x,y,t)=v(x-ct),\] 
where
$c\in\R$ is the speed of propagation, and $v$ satisfies the ODE
\[
\left(c v^\prime+v^{\prime\prime\prime}+vv^\prime\right)^\prime=0.
\]
Integrating this equation twice, and writing $x$ instead of $x-ct$, we
obtain the second order ODE
\[
v^{\prime\prime}= -cv -\frac12v^2+b+dx,
\]
in which $b$ and $d$ are arbitrary constants. Since we are interested
in periodic solutions, we can set $d=0$ and the equation becomes
\begin{equation}\label{e:onedtw}
v^{\prime\prime}= -cv -\frac12v^2+b.
\end{equation}

\paragraph{Symmetries} The  KP equation \eqref{e:kp} possesses the
\emph{Galilean invariance}
\[
u(x,y,t)\mapsto u(x+\alpha t,y,t) +\alpha,
\]
which leads to the invariance of the traveling-wave equation
\eqref{e:onedtw} under the transformation
\[
v\mapsto v+\alpha,\quad c\to c-\alpha,\quad b\mapsto b+\alpha
c-\frac12\alpha^2. 
\]
As a consequence, we can restrict to positive speeds $c>0$, and to 
$b=0$. 
In addition, the
KP equation \eqref{e:kp} possesses the 
\emph{scaling invariance}
\[
u(x,y,t)\mapsto c u(\sqrt c x, cy, c\sqrt c t), 
\]
which allows to restrict to the speed $c=1$.

\paragraph{Small periodic waves} Restricting to $c=1$ and $b=0$, the
equation \eqref{e:onedtw} becomes
\begin{equation}\label{e:onedtw0}
v^{\prime\prime}= -v -\frac12v^2.
\end{equation}
In the phase-plane $(v,v^\prime)$, this equation possesses the
equilibria $(0,0)$ and $(-2,0)$, which are a center and a saddle,
respectively. The center $(0,0)$ is surrounded by a one-parameter
family of periodic orbits, with periods $L\in(2\pi,\infty)$.  As $L\to
2\pi$, the periodic orbits shrink to the equilibrium $(0,0)$, whereas
for $L\to\infty$ the periodic orbits tend to a homoclinic orbit
connecting the saddle equilibrium $(-2,0)$ to itself.  These periodic
orbits are symmetric with respect to the $v$-axis in the  phase-plane
$(v,v^\prime)$, so that to each periodic orbit corresponds a 
periodic solution of 
\eqref{e:onedtw0} which is even in $x$, and unique up to translations
in $x$.

We are
interested here in the periodic orbits close to the center equilibrium
$(0,0)$, which correspond to small periodic traveling waves of the KP
equation \eqref{e:kp}. Looking for solutions of 
\eqref{e:onedtw0} of the form $v(x)=P(kx)$, where $k$ is the
wavenumber and $P$ is a $2\pi$-periodic function, a direct calculation
shows that the one-parameter family of small periodic solutions of 
\eqref{e:onedtw0}, which are even in $x$, is given by
\[
v_a(x)=P_a(k_a x),
\]
where
\begin{equation}\label{e:expptw}
P_a(z) = a\cos(z)+\frac14\left(\frac13\cos(2z)-1\right)a^2+O(|a|^3),
\quad k_a^2=1-\frac5{24}a^2+O(a^4).
\end{equation}
The parameter $a$ is small and corresponds to the Fourier mode
$1$ in the Fourier expansion of the $2\pi$-periodic function $P_a$,
\[
a=\frac1\pi\int_0^{2\pi} P_a(z)\cos(z)\,dz.
\]
The function $P_a$ is even in $z$, and notice the symmetry
\[
P_a(z+\pi)=P_{-a}(z),\quad k_a = k_{-a}.
\]
Taking into account the translation invariance we obtain a
two-parameter family of periodic solutions of \eqref{e:onedtw0}, and
including the parameters $b$ and $c$, we find a four parameter
family of one-dimensional periodic traveling waves of the KP equation
\eqref{e:kp}.

\section{The spectral stability problem}
\label{s:3}

Consider the one-dimensional periodic traveling waves of the KP
equation \eqref{e:kp} found in Section~\ref{s:2}, 
\[
u_a(x,y,t)=P_a(k_a(x-t)).
\]
We introduce the scaling
\[
z=k_a(x-t),\quad \tilde y=k_a^2y,\quad \tilde t=k_a^3t,
\]
which transforms the KP equation \eqref{e:kp} into
\begin{equation}\label{e:kps}
u_{tz}-u_{zzzz}-\frac1{k_a^2}u_{zz}-\frac1{k_a^2}(uu_z)_z+\sigma u_{yy}=0,
\end{equation}
where we have dropped the tilde for notational simplicity. Then $P_a$
is a stationary solution of \eqref{e:kps}, and our goal is to study
its spectral stability. 

Following the standard approach to spectral stability, we linearize 
\eqref{e:kps} about the stationary
solution $P_a$ and obtain the linear evolution equation
\begin{equation}\label{e:kplin}
w_{tz}-w_{zzzz}-\frac1{k_a^2}w_{zz}-\frac1{k_a^2}(P_a
w)_{zz}+\sigma w_{yy}=0.
\end{equation}
This equation has coefficients depending upon $z$, only, so that we make the
Ansatz,
\[
w(z,y,t)=e^{\lambda t+i\ell y}W(z),
\]
which leads to the equation
\[
\lambda W_{z}-W_{zzzz}-\frac1{k_a^2}W_{zz}-\frac1{k_a^2}(P_a
W)_{zz}-\sigma\ell^2 W=0.
\]
The left hand side of this equation defines the linear differential operator
\[
\mathcal M_a(\lambda,\ell) = \lambda\partial_z-\partial_{z}^4-
\frac1{k_a^2}\partial_{z}^2((1+P_a)\cdot)-\sigma\ell^2,
\]
and the spectral stability problem is concerned with the invertibility
of this operator for $\lambda\in\C$ and $\ell\in\R$: the periodic wave
is spectrally stable if this operator is invertible for any
$\lambda\in\C$ with $\Re\lambda>0$, and unstable otherwise. The type of the
allowed perturbations is determined by the choice of the function
space and the values of $\ell$.

\paragraph{One-dimensional perturbations}
The spectral stability problem for 
one-dimensional perturbations, i.e., perturbations
which are independent of $y$, corresponds to the case
$\ell=0$, when
\[
\mathcal M_a(\lambda,0) = \partial_z\mathcal K_a(\lambda),\quad 
\mathcal K_a(\lambda)=
\lambda-\partial_{z}^3-
\frac1{k_a^2}\partial_{z}((1+P_a)\cdot).
\]
The linear operator $\mathcal K_a(\lambda)$ in this decomposition is
in fact the linear operator arising in the spectral stability problem
for the KdV equation \eqref{e:kdv}. As solutions of the KdV equation,
the periodic waves $P_a$ are spectrally stable with respect to
perturbations which are periodic in $z$, but also with respect to 
perturbations 
which are localized or bounded in $z$, i.e., the linear
operator $\mathcal K_a(\lambda)$ is invertible, for any
$\lambda\in\C$, $\Re\lambda>0$, in $L^2(0,2\pi)$, and also $L^2(\R)$
or $C_b(\R)$ (e.g., see \cite{HK,BD}). 
On the other hand, the operator $\partial_z$ is not invertible in
these spaces, so that the linear operator $\mathcal M_a(\lambda,0)$ is
not invertible, for any $\lambda\in\C$, $\Re\lambda>0$. This
ill-posedness of the spectral stability problem for one-dimensional
perturbations shows that for such perturbations
we have to use the operator $\mathcal K_a(\lambda)$ instead of
$\mathcal M_a(\lambda,0)$, i.e., the KP operator $\mathcal
M_a(\lambda,\ell)$ should be replaced by the KdV operator $\mathcal
K_a(\lambda)$, when $\ell=0$.

\paragraph{Two-dimensional perturbations}
Truly two-dimensional perturbations correspond to values
$\ell\not=0$. We consider three types of such perturbations:
perturbations which are periodic in $z$, perturbations which are
localized in $z$, and perturbations which are bounded in $z$. The type
of these perturbations is determined by choice of the function space
in which acts the linear operator $\mathcal M_a(\lambda,\ell)$:
\begin{enumerate}
\item  $\mathcal M_a(\lambda,\ell)$ is considered in  $L^2(0,2\pi)$,
  with domain
\[
H^4_\per(0,2\pi) =\{f\in H^4_\loc(\R)\;;\; f(z)=f(z+2\pi),\
\forall\ z\in\R\},
\]
for perturbations which are periodic in $z$;
\item  $\mathcal M_a(\lambda,\ell)$ is considered in  $L^2(\R)$,
  with domain $H^4(\R)$,
for perturbations which are localized in $z$;
\item  $\mathcal M_a(\lambda,\ell)$ is considered in  $C_b(\R)$,
  with domain $C_b^4(\R)$,
for perturbations which are bounded in $z$.
\end{enumerate}

\paragraph{Spectral stability} Summarizing, we give the following
  definition.

\begin{Definition}[Spectral stability]
\begin{enumerate}
\item 
We say that the periodic wave $P_a$ is \emph{spectrally stable in one
  dimension}  
\emph{with respect to periodic perturbations} (resp. \emph{localized
  or bounded perturbations}), if the KdV operator
  $\mathcal K_a(\lambda)$ acting in $L^2(0,2\pi)$ (resp.  $L^2(\R)$ or
  $C_b(\R)$) 
  with domain
  $H^3_\per(0,2\pi)$ (resp.  $H^3(\R)$ or   $C_b^3(\R)$) is
  invertible, for any $\lambda\in\C$, $\Re\lambda>0$. 

\item
We say that the periodic wave $P_a$ is \emph{transversely spectrally
  stable} \emph{with respect to periodic perturbations} (resp. \emph{localized
  or bounded perturbations}), if  it is spectrally stable in one
  dimension, and the KP operator
  $\mathcal M_a(\lambda,\ell)$ acting in $L^2(0,2\pi)$ (resp.  $L^2(\R)$ or
  $C_b(\R)$)  with domain
  $H^4_\per(0,2\pi)$  (resp.  $H^4(\R)$ or   $C_b^4(\R)$)
is invertible, for any $\lambda\in\C$, $\Re\lambda>0$ and any $\ell\not=0$;
\end{enumerate}
\end{Definition}

The results in \cite{HK,BD} show that the periodic waves $P_a$ are
spectrally stable in one dimension. We discuss their transverse
spectral stability in the following sections.  It turns out that the
stability properties strongly depend upon the sign of $\sigma$. We
focus on the case $\sigma=1$, i.e., the KP-I equation, in
Sections~\ref{s:4} and \ref{s:5}, and then briefly discuss 
the case  $\sigma=-1$, i.e., the
KP-II equation, in Section~\ref{s:6}. 
Notice that depending upon the values of $\ell$ we may distinguish three
different regimes: {\it short wavelength transverse perturbations}, when
$\ell\gg1$, {\it long wavelength transverse perturbations}, when
$\ell\ll1$, and {\it finite wavelength transverse perturbations}, otherwise.
The results in \cite{JZ} concern the transverse spectral stability with
respect to periodic perturbations in the regime of long wavelength
transverse perturbations.

\section{KP-I equation: periodic perturbations}
\label{s:4}

In this section we study the transverse spectral stability of the
periodic waves $P_a$ with respect to periodic perturbations, for the
KP-I equation, i.e., when $\sigma=1$. More precisely, we
study the invertibility of the operator 
\[
\mathcal M_a(\lambda,\ell) = \lambda\partial_z-\partial_{z}^4-
\frac1{k_a^2}\partial_{z}^2((1+P_a)\cdot)-\ell^2,
\]
acting in $L^2(0,2\pi)$ with domain $H^4_\per(0,2\pi)$,
for $\lambda\in\C$, $\Re\lambda>0$, and $\ell\in\R$, $\ell\not=0$.

\subsection{Reformulation of the spectral stability problem}
\label{ss:specper}

We show that the question of the invertibility of the operator
$\mathcal M_a(\lambda,\ell)$ is equivalent to the study of the
spectrum of the linear operator
\[
\mathcal A_a(\ell) = \partial_{z}^3+
\frac1{k_a^2}\partial_{z}((1+P_a)\cdot)+\ell^2
\partial_z^{-1}, 
\]
acting in the space
\[
L^2_0(0,2\pi)=\left\{f\in L^2(0,2\pi)\;;\;\int_0^{2\pi}f(z)\,dz=0\right\},
\]
of square-integrable functions on $(0,2\pi)$ with of zero-mean, with
domain $H^3_\per(0,2\pi)\cap L^2_0(0,2\pi)$. Here
$\partial_z^{-1}$ is the inverse of the restriction of $\partial_z$ to
the subspace $L^2_0(0,2\pi)$.

\begin{Lemma}
Assume that $\lambda\in\C$ and  $\ell\in\R$, $\ell\not=0$. Then the linear
operator $\mathcal M_a(\lambda,\ell)$ 
acting in $L^2(0,2\pi)$ with domain $H^4_\per(0,2\pi)$ is invertible
if and only if its restriction to the subspace
$L^2_0(0,2\pi)$ is an invertible operator.
\end{Lemma}

\begin{Proof}
First, notice that the subspace $L^2_0(0,2\pi)\subset L^2(0,2\pi)$ is
invariant under the action of $\mathcal M_a(\lambda,\ell)$, since 
\[
\int_0^{2\pi}(\mathcal M_a(\lambda,\ell)W)(z)\,dz = -\ell^2 
\int_0^{2\pi} W(z)\,dz=0,
\]
for any $W\in H^4_\per(0,2\pi)\cap L^2_0(0,2\pi)$. Next, the
operator  $\mathcal M_a(\lambda,\ell)$ has compact resolvent, since
$H^4_\per(0,2\pi)$  is compactly embedded in
$L^2(0,2\pi)$. Consequently, the spectrum of $\mathcal
M_a(\lambda,\ell)$ consists of isolated eigenvalues with finite
algebraic multiplicity. In particular, $\mathcal
M_a(\lambda,\ell)$  is
invertible if and only if $0$ is an eigenvalue of  $\mathcal
M_a(\lambda,\ell)$, i.e., if and only if there exists $W\in H^4_\per(0,2\pi)$,
$W\not=0$, such that 
\[
\mathcal M_a(\lambda,\ell)W=0.
\]
Since $\ell\not=0$, we conclude that $W\in L^2_0(0,2\pi)$, so that any
eigenfunction $W\in H^4_\per(0,2\pi)$ associated to the eigenvalue $0$
belongs to $L^2_0(0,2\pi)$. This implies that $0$ is an eigenvalue of
$\mathcal M_a(\lambda,\ell)$ if and only if $0$ is an eigenvalue of
the restriction of  $\mathcal M_a(\lambda,\ell)$ to $L^2_0(0,2\pi)$,
which proves the lemma.
\end{Proof}

Since the operator $\partial_z$ acting in  $L^2_0(0,2\pi)$ with domain
$H^1_\per(0,2\pi)\cap L^2_0(0,2\pi)$ is invertible, the
following result is an immediate consequence of the above lemma.

\begin{Corollary}\label{c:inv}
Assume that $\lambda\in\C$ and  $\ell\in\R$, $\ell\not=0$. Then the linear
operator $\mathcal M_a(\lambda,\ell)$ 
acting in $L^2(0,2\pi)$ with domain $H^4_\per(0,2\pi)$ is invertible
if and only if $\lambda$ belongs to the spectrum of the operator
$\mathcal A_a(\ell)$ acting in $L^2_0(0,2\pi)$ with domain
$H^3_\per(0,2\pi)\cap L^2_0(0,2\pi)$. 
\end{Corollary}

Our problem consists now in the study of the spectrum of
the operator $\mathcal A_a(\ell)$ acting in $L^2_0(0,2\pi)$ with domain
$H^3_\per(0,2\pi)\cap L^2_0(0,2\pi)$, for $\ell\in\R$, $\ell\not=0$.  
This operator has compact resolvent, just as  $\mathcal
M_a(\lambda,\ell)$, so that its spectrum consists of isolated
eigenvalues with finite algebraic multiplicity. 
In addition, the spectrum has the following symmetry property.

\begin{Lemma}\label{l:syms}
The spectrum of  $\mathcal
A_a(\ell)$ is symmetric with respect to both the real and the
imaginary axis.
\end{Lemma}

\begin{Proof}
First, the spectrum of  $\mathcal A_a(\ell)$ is symmetric with respect
to the real axis, since $\mathcal A_a(\ell)$ is a real operator. Next,
consider  the reflection
$\mathcal S$ defined by
\begin{equation}\label{e:s}
\mathcal SW(z) = W(-z),
\end{equation}
and notice that $\mathcal A_a(\ell)$ anti-commutes with
$\mathcal S$,
\[
(\mathcal A_a(\ell) \mathcal S W)(z) = \mathcal A_a(\ell)
(W(-z)) = -(\mathcal A_a(\ell)W)(-z)= -(\mathcal S\mathcal A_a(\ell)W)(z),
\]
where we have used the fact that $P_a$ is an
even function. If $\lambda$ is an eigenvalue of  $\mathcal A_a(\ell)$ with
associated eigenvector $W_\lambda$, 
\[
\mathcal A_a(\ell) W_\lambda=\lambda W_\lambda,
\]
then 
\[
\mathcal A_a(\ell) \mathcal SW_\lambda=
-\mathcal S\mathcal A_a(\ell)W_\lambda = 
-\lambda \mathcal SW_\lambda.
\]
Consequently, $-\lambda$ is an eigenvalue of $\mathcal
A_a(\ell)$. This implies that the spectrum of $\mathcal A_a(\ell)$ is also
symmetric with respect to the origin, and completes the proof.
\end{Proof}

The analysis of the spectrum of $\mathcal A_a(\ell)$ in the next
sections relies upon two main tools:
\begin{itemize}
\item
the decomposition of  $\mathcal
A_a(\ell)$ into the product of a skew-adjoint and a self-adjoint
operator:
\begin{equation}\label{e:dec}
\mathcal A_a(\ell) = -\partial_z\mathcal L_a(\ell),\quad 
\mathcal L_a(\ell)=-\partial_{z}^2-
\frac1{k_a^2}\left((1+P_a)\cdot\right)-\ell^2\partial_z^{-2}, 
\end{equation}
which  leads to a simple characterization of the unstable eigenvalues of
$\mathcal A_a(\ell)$; 
\item
perturbation arguments for linear operators: we regard  $\mathcal
A_a(\ell)$ as a perturbation, for small $a$, 
of the operator with constant
coefficients 
\[
\mathcal A_0(\ell) = \partial_{z}^3+
\partial_{z}+\ell^2
\partial_z^{-1}, 
\]
obtained by setting $a=0$.
\end{itemize}
The main result is the following theorem showing the
spectral instability of the periodic waves $P_a$, for $a$ sufficiently
small. The instability is due to a pair of real eigenvalues, with
opposite signs, which is found in the long wavelength
regime, for $\ell^2=O(a^2)$.

\begin{Theorem}\label{t:per} For any $a$ sufficiently small, 
there exists $\ell_a^2 =
\frac1{12}a^2+O(a^4)$, such that 
\begin{enumerate}
\item for any $\ell^2\geq\ell_a^2$,
the spectrum of $\mathcal
A_a(\ell)$ is purely imaginary;
\item for any $\ell^2<\ell_a^2$,
the spectrum of $\mathcal
A_a(\ell)$ is purely imaginary, except for a pair of  simple real
eigenvalues, with opposite signs.
\end{enumerate}
\end{Theorem}

\subsection{Characterization of the unstable spectrum}

Consider the decomposition \eqref{e:dec} of the linear operator
$\mathcal A_a(\ell)$. 
The operators $-\partial_z$ and  $\mathcal L_a(\ell)$ in this
decomposition are skew-adjoint and self-adjoint, respectively, as
operators acting in $L^2(0,2\pi)$. 
In  \cite{HK} it has been shown, under some general assumptions, that
for an operator $\mathcal A=\mathcal
J\mathcal L$, where $\mathcal J$ and $\mathcal L$ are skew-adjoint and
self-adjoint, respectively,
the number of unstable eigenvalues, counted
with multiplicities, is not larger than the number of negative
eigenvalues of the operator 
$\mathcal L$. In particular, if $\mathcal L$ has no
negative eigenvalues, then $\mathcal A$ has no unstable
eigenvalues. 
For the operators in \eqref{e:dec} we cannot directly apply
this result, since the subspace $L^2_0(0,2\pi)$ on which
acts $\mathcal A_a(\ell)$ is not an invariant subspace for the 
operator  $\mathcal L_a(\ell)$. Nevertheless, we can use part of the
arguments in \cite{HK} to prove the following result.

\begin{Lemma}\label{l:uev}
Assume that $\lambda$ is an eigenvalue of the operator $\mathcal
A_a(\ell)$ acting in $L^2_0(0,2\pi)$, and
that $W_\lambda$ is an associated eigenvector. If $\Re\lambda\not=0$,
then
\[
\langle\mathcal L_a(\ell)W_\lambda,W_\lambda\rangle = 0,
\]
where $\langle\cdot,\cdot\rangle$ is the usual scalar product in
$L^2(0,2\pi)$.
\end{Lemma}

\begin{Proof} Recall that the operator $\partial_z$ acting on
  $L^2_0(0,2\pi)$ is skew-adjoint and invertible.
Then, following the proof of \cite[Lemma 2.7]{HK}, we find
\[
\langle\mathcal L_a(\ell)W_\lambda,W_\lambda\rangle =
\langle \mathcal L_a(\ell)W_\lambda,
\partial_z\partial_z^{-1}W_\lambda\rangle =
\langle\mathcal A_a(\ell)W_\lambda,\partial_z^{-1}W_\lambda\rangle =
\lambda\langle W_\lambda,\partial_z^{-1}W_\lambda\rangle, 
\]
and similarly,
\[
\langle W_\lambda,\mathcal L_a(\ell)W_\lambda\rangle =
\langle\partial_z \partial_z^{-1}W_\lambda,
\mathcal L_a(\ell)W_\lambda\rangle =
\langle\partial_z^{-1}W_\lambda,
\mathcal A_a(\ell)W_\lambda\rangle =
\overline\lambda\langle \partial_z^{-1}W_\lambda,W_\lambda\rangle
=-\overline\lambda\langle W_\lambda,\partial_z^{-1}W_\lambda\rangle. 
\]
Since $\mathcal L_a(\ell)$ is self-adjoint, we conclude that
\[
(\lambda+\overline\lambda)\langle\mathcal
L_a(\ell)W_\lambda,W_\lambda\rangle = 0,
\]
so that $\langle\mathcal
L_a(\ell)W_\lambda,W_\lambda\rangle = 0$, since $\Re\lambda\not=0$.
\end{Proof}

\begin{Corollary}\label{c:imspec}
Assume that there exists a positive constant $c$ such that
\begin{equation}\label{e:posL}
\langle\mathcal
L_a(\ell)W,W\rangle \geq c\langle W,W\rangle,\quad\forall\ W\in
H^3_\per(0,2\pi)\cap L^2_0(0,2\pi).
\end{equation}
Then the spectrum of $\mathcal A_a(\ell)$ is purely imaginary.
\end{Corollary}

\begin{Proof}
Assuming that there exists an eigenvalue $\lambda$ of $\mathcal
A_a(\ell)$ with $\Re\lambda\not=0$, then any
associated eigenvector $W_\lambda\not=0$
satisfies $\langle\mathcal
L_a(\ell)W_\lambda,W_\lambda\rangle >0$, by hypothesis. This
contradicts the result in  Lemma~\ref{l:uev}, and proves the corollary.
\end{Proof}

\subsection{Finite and short wavelength transverse perturbations}

We start the analysis of the spectrum of $\mathcal A_a(\ell)$ with the
values of $\ell$ away from the origin, $|\ell|\geq\ell_*$, for some
$\ell_*>0$, i.e., finite and  short wavelength transverse perturbations.
We use the result in Corollary~\ref{c:imspec} to show that  the spectrum of
$\mathcal A_a(\ell)$ is purely imaginary for such values of $\ell$,
provided $a$ is sufficiently small.

\begin{Lemma}\label{l:short}
Assume that $\ell_*>0$. There exists $a_*>0$, such that the  spectrum of
$\mathcal A_a(\ell)$ is purely imaginary, for any $\ell$ and $a$
satisfying $|\ell|\geq\ell_*$ and $|a|\leq a_*$.
\end{Lemma}

\begin{Proof}
According to the Corollary~\ref{c:imspec} it is enough to show that
\begin{equation}\label{e:cstar}
\langle\mathcal
L_a(\ell)W,W\rangle \geq c_*\langle W,W\rangle,\quad\forall\ W\in
H^3_\per(0,2\pi)\cap L^2_0(0,2\pi),
\end{equation}
for some $c_*>0$.
We write
\[
\mathcal L_a(\ell) = \mathcal L_0(\ell) + \widetilde{\mathcal L}_a,\quad
\mathcal L_0(\ell)=-\partial_{z}^2-1-\ell^2\partial_z^{-2},\quad
\widetilde{\mathcal L}_a =1-
\frac1{k_a^2}\left(1+P_a\right).
\]
The operator $\mathcal L_0(\ell)$ has constant coefficients, and using
Fourier series we find that its spectrum in $L^2_0(0,2\pi)$ is given
by
\[
\sigma(\mathcal L_0(\ell))=
\left\{\mu_n(\ell)=n^2-1+\frac{\ell^2}{n^2}\;;\;n\in\Z^*\right\},
\]
and satisfies
\[
\sigma(\mathcal L_0(\ell))\subset[c_0,\infty),\quad c_0=\min\{\ell_*^2,1\},
\]
for any $\ell$, $|\ell|\geq\ell_*$. As a consequence, we have that
\begin{equation}\label{e:c0}
\langle\mathcal
L_0(\ell)W,W\rangle \geq c_0\langle W,W\rangle,\quad\forall\ W\in
H^3_\per(0,2\pi)\cap L^2_0(0,2\pi).
\end{equation}
The linear operator $\widetilde{\mathcal L}_a$ is bounded in
$L^2(0,2\pi)$, and a direct calculation shows that 
\[
\|\widetilde{\mathcal L}_a\|\leq c_1 a,  
\]
for some $c_1>0$, and any $a$ sufficiently small. Together with
\eqref{e:c0}, this implies that
\[
\langle\mathcal
L_a(\ell)W,W\rangle = 
\langle\mathcal
L_0(\ell)W,W\rangle + \langle\widetilde{\mathcal L}_a W,W\rangle 
\geq (c_0-c_1a)\langle W,W\rangle,\quad\forall\ W\in
H^3_\per(0,2\pi)\cap L^2_0(0,2\pi),
\]
which proves \eqref{e:cstar}, for any $\ell$, $|\ell|\geq\ell_*$, 
provided $a$ is sufficiently small.
\end{Proof}

\subsection{Long wavelength transverse perturbations}\label{ss:longper}

We consider now  the spectrum of $\mathcal A_a(\ell)$ for small
$\ell$. Notice that the arguments in the proof of 
Lemma~\ref{l:short} do not work when $\ell$ is small, since the
constant $c_0$ in \eqref{e:c0} tends to zero, as $\ell\to0$. 

Here, we regard  $\mathcal A_a(\ell)$ as a perturbation of the
operator with constant coefficients
\[
\mathcal A_0(0) = \partial_z^3+\partial_z,
\]
acting in $L^2_0(0,2\pi)$. 
A direct calculation shows that the spectrum of $\mathcal A_0(0)$ is
given by
\[
\sigma(\mathcal A_0(0))=\{i\omega_n=-in^3+in\;;\;n\in\Z^*\}.
\]
In particular, zero is a double eigenvalue of   $\mathcal A_0(0)$, and
the remaining eigenvalues are all simple, purely imaginary, and
located outside the open ball $B(0;5)$  of radius $5$ centered at the origin.

\begin{Lemma}\label{l:spdec}
The following properties hold, for any $\ell$ and $a$ sufficiently
small.
\begin{enumerate}
\item The spectrum of $\mathcal A_a(\ell)$ decomposes as
\[
\sigma(\mathcal A_a(\ell)) = \sigma_0(\mathcal A_a(\ell))\cup 
\sigma_1(\mathcal A_a(\ell)),
\]
with
\[
\sigma_0(\mathcal A_a(\ell)) \subset B(0;1),\quad 
\sigma_1(\mathcal A_a(\ell)) \subset \C\setminus \overline{B(0;4)}. 
\]
\item The spectral projection $\Pi_a(\ell)$ associated with 
$\sigma_0(\mathcal A_a(\ell)) $ satisfies
  $\|\Pi_a(\ell)-\Pi_0(0)\|=O(\ell^2+|a|)$.
\item The spectral subspace $\mathcal X_a(\ell) =
  \Pi_a(\ell)(L^2_0(0,2\pi))$ is two dimensional.
\end{enumerate}
\end{Lemma}

\begin{Proof}
{\it (i)} Consider $\lambda\in \overline{B(0;4)}\setminus
B(0;1)$. Then $\lambda$ belongs to the resolvent set of $\mathcal
A_0(0)$, i.e., the operator $ \lambda -\mathcal A_0(0)$ is
invertible. We write
\[
\lambda -\mathcal A_a(\ell) = \left( I - 
\widetilde{\mathcal A}_a(\ell) (\lambda -\mathcal A_0(0))^{-1} \right) 
(\lambda-\mathcal A_0(0)), 
\]
where $\widetilde{\mathcal A}_a(\ell) =  \mathcal A_a(\ell) - \mathcal
A_0(0)$.
A straightforward calculation shows that there exist positive
constants $c_0$ and $c_1$ such that
\[
\|(\lambda -\mathcal A_0(0))^{-1}\|_{L^2\to H^1}\leq c_0,\quad
\|\widetilde{\mathcal A}_a(\ell)\|_{H^1\to L^2}\leq c_1(\ell^2+|a|),
\]
for any  $\lambda\in \overline{B(0;4)}\setminus
B(0;1)$, and any $\ell$ and $a$ sufficiently small. Consequently,
\[
\|\widetilde{\mathcal A}_a(\ell)(\lambda -\mathcal
A_0(0))^{-1}\|\leq c_0  c_1(\ell^2+|a|)\leq \frac12,
\]
for any  $\lambda\in \overline{B(0;4)}\setminus
B(0;1)$, provided $\ell$ and $a$ are sufficiently small, so that the
operator $ I - 
\widetilde{\mathcal A}_a(\ell) (\lambda -\mathcal A_0(0))^{-1}$ is
invertible. Consequently, $\lambda -\mathcal A_a(\ell)$ is invertible,
for any $\lambda\in \overline{B(0;4)}\setminus B(0;1)$, which proves the first
part of the lemma.

{\it (ii)} The spectral projection $\Pi_a(\ell)$ can be computed with
the help of the Dunford integral formula
\[
\Pi_a(\ell) = \frac1{2\pi i}\int_{\partial B(0;1)} 
\left( \lambda -\mathcal A_a(\ell)\right)^{-1}\,d\lambda,
\]
where, according to the previous arguments,
\begin{eqnarray*}
\left( \lambda -\mathcal A_a(\ell)\right)^{-1} &=& 
\left( \lambda -\mathcal A_0(0)\right)^{-1}
\left( I - 
\widetilde{\mathcal A}_a(\ell) (\lambda -\mathcal A_0(0))^{-1}
\right)^{-1}\\
&=& \left( \lambda -\mathcal A_0(0)\right)^{-1} 
\sum_{k\geq0} \left(\widetilde{\mathcal A}_a(\ell) (\lambda -\mathcal
A_0(0))^{-1}\right)^k\\
&=& \left( \lambda -\mathcal A_0(0)\right)^{-1} +  \left( \lambda
-\mathcal A_0(0)\right)^{-1} 
\sum_{k\geq1} \left(\widetilde{\mathcal A}_a(\ell) (\lambda -\mathcal
A_0(0))^{-1}\right)^k.
\end{eqnarray*}
Consequently,
\[
\Pi_a(\ell) - \Pi_0(0)=  \frac1{2\pi i}\int_{\partial B(0;1)} 
 \left( \lambda
-\mathcal A_0(0)\right)^{-1} 
\sum_{k\geq1} \left(\widetilde{\mathcal A}_a(\ell) (\lambda -\mathcal
A_0(0))^{-1}\right)^k\,d\lambda,
\]
and since
\[
\big\|
\left( \lambda
-\mathcal A_0(0)\right)^{-1} 
\sum_{k\geq1} \left(\widetilde{\mathcal A}_a(\ell) (\lambda -\mathcal
A_0(0))^{-1}\right)^k
\big\|
\leq 2 c_0^2  c_1(\ell^2+|a|),
\]
for any $\lambda\in\partial B(0;1)$, and any $\ell$ and $a$
sufficiently small, we conclude that  
$\|\Pi_a(\ell)-\Pi_0(0)\|=O(\ell^2+|a|)$.

{\it (iii)} The last part of the lemma is an immediate consequence of
the estimate in (ii) (e.g., see \cite[Lemma B.1]{HLS}).
\end{Proof}

Next, we use  the
result in Lemma~\ref{l:uev}, 
to locate the eigenvalues in $\sigma_1(\mathcal A_a(\ell))$, 
and then we determine the two
eigenvalues in 
$\sigma_0(\mathcal A_a(\ell))$, by computing an expansion of these
eigenvalues for small $\ell$ and $a$. These results complete the proof
of Theorem~\ref{t:per}.

\begin{Lemma}\label{l:long1}
Consider the decomposition of the spectrum of $\mathcal A_a(\ell)$ in
Lemma~\ref{l:spdec}.  Then
$\sigma_1(\mathcal A_a(\ell))\subset i\R$, for any  $\ell$ and $a$ 
sufficiently small. 
\end{Lemma}

\begin{Proof}
Consider the restriction $\mathcal A_a^1(\ell)$ of $\mathcal
A_a(\ell)$ to the spectral subspace $\mathcal Y_a(\ell) =
(I-  \Pi_a(\ell))(L^2_0(0,2\pi))$, so that
\[
\sigma_1(\mathcal A_a(\ell))= \sigma(\mathcal A_a^1(\ell)).
\]
Assume that $\lambda$ is an eigenvalue of  $\mathcal A_a^1(\ell)$, and
that $W_\lambda\not=0$ is an associated eigenvector. In particular,
$W_\lambda = (I-\Pi_a(\ell))W_\lambda$. We use the result
in  Lemma~\ref{l:uev} to show that $\Re\lambda=0$.

Recall that $\mathcal A_a(\ell)=-\partial_z\mathcal L_a(\ell)$, and
write
\[
\mathcal L_a(\ell)=\mathcal L_0(0)+\widetilde{\mathcal L}_a(\ell),\quad
\mathcal L_0(0)=-\partial_{z}^2-1,\quad
\widetilde{\mathcal L}_a(\ell) =1-
\frac1{k_a^2}\left(1+P_a\right)-\ell^2\partial_z^{-2},
\]
so that $\widetilde{\mathcal L}_a(\ell)$ is a bounded operator on
$L^2(0,2\pi)$, with
norm
\begin{equation}\label{e:estl}
\|\widetilde{\mathcal L}_a(\ell)\|\leq c_0(\ell^2+|a|),
\end{equation}
for some positive constant $c_0$, and any $\ell$ and $a$ sufficiently
small. We set $\widetilde\Pi_a(\ell)=\Pi_a(\ell)-\Pi_0(0)$, and then
the property in Lemma~\ref{l:spdec}~(ii) shows that
\begin{equation}\label{e:estpi}
\|\widetilde\Pi_a(\ell)\|=
\|\Pi_a(\ell)-\Pi_0(0)\|\leq c_1(\ell^2+|a|),
\end{equation}
for some positive constant $c_1$, and any $\ell$ and $a$ sufficiently
small. We compute
\begin{eqnarray*}
\langle\mathcal L_a(\ell) W_\lambda,W_\lambda\rangle&=&
\langle\mathcal L_a(\ell) (I-\Pi_a(\ell))
W_\lambda,(I-\Pi_a(\ell))W_\lambda\rangle \\
&=& \langle\mathcal L_0(0) (I-\Pi_a(\ell))
W_\lambda,(I-\Pi_a(\ell))W_\lambda\rangle +
\langle\widetilde{\mathcal L}_a(\ell) (I-\Pi_a(\ell))
W_\lambda,(I-\Pi_a(\ell))W_\lambda\rangle   \\
&=& \langle\mathcal L_0(0) (I-\Pi_0(0))
W_\lambda,(I-\Pi_0(0))W_\lambda\rangle -
\langle\mathcal L_0(0) \widetilde\Pi_a(\ell)
W_\lambda,(I-\Pi_0(0))W_\lambda\rangle \\
&&-\langle\mathcal L_0(0) (I-\Pi_a(\ell))
W_\lambda,\widetilde\Pi_a(\ell)W_\lambda\rangle +
\langle\widetilde{\mathcal L}_a(\ell) (I-\Pi_a(\ell))
W_\lambda,(I-\Pi_a(\ell))W_\lambda\rangle  .
\end{eqnarray*}
Since the spectrum of the restriction of $\mathcal L_0(0)$ to
$(I-\Pi_0(0))L^2_0(0,2\pi)$ consists of the eigenvalues
\[
n^2-1,\quad n\in\Z\setminus\{-1,0,1\},
\]
we have that
\[
\langle\mathcal L_0(0) (I-\Pi_0(0))
W_\lambda,(I-\Pi_0(0))W_\lambda\rangle \geq 3\langle
W_\lambda,W_\lambda\rangle, 
\]
and taking into account the estimates \eqref{e:estl} and
\eqref{e:estpi}, we conclude that 
\[
\langle\mathcal L_a(\ell) W_\lambda,W_\lambda\rangle\geq
(3-c_2(\ell^2+|a|)) \langle
W_\lambda,W_\lambda\rangle, 
\]
for some positive constant $c_2$, and any $\ell$ and $a$ sufficiently
small. Consequently, $\langle\mathcal L_a(\ell)
W_\lambda,W_\lambda\rangle>0$, provided  $\ell$ and $a$ are sufficiently
small, so that $\Re\lambda=0$, by Lemma~\ref{l:uev}.
\end{Proof}

\begin{Lemma}\label{l:long0}
Assume that $\ell$ and $a$ are sufficiently small. There exists $\ell_a^2 =
\frac1{12}a^2+O(a^4)$, such that the two eigenvalues in $\sigma_0(\mathcal
A_a(\ell))$ are purely imaginary, if $\ell^2\geq\ell_a^2$, and are real,
one positive and one negative, if $\ell^2<\ell_a^2$.
\end{Lemma}

\begin{Proof}
The eigenvalues in
$\sigma_0((\mathcal A_a(\ell))$ are the eigenvalues of 
the restriction of $\mathcal
A_a(\ell)$ to the two-dimensional spectral subspace $\mathcal X_a(\ell) =
\Pi_a(\ell)(L^2_0(0,2\pi))$. We determine the location of these
eigenvalues by computing successively
a basis of  $\mathcal X_a(\ell)$, 
the $2\times2$ matrix representing the action of $\mathcal A_a(\ell)$
on this basis, and the eigenvalues of this matrix.

We start with the computation for $a=0$. Then 
\[
\mathcal A_0(\ell)= \partial_{z}^3+
\partial_{z}+\ell^2
\partial_z^{-1}, 
\] 
is an operator with constant coefficients, and 
\[
\sigma_0(\mathcal A_0(\ell))=\{-i\ell^2,i\ell^2\}.
\]
The associated
eigenvectors are $e^{iz}$ and $e^{-iz}$, and we choose 
\[
\xi_0^0(\ell) = \sin(z),\quad \xi_0^1(\ell)  =\cos(z),
\]
as  basis of the corresponding spectral subspace. Since 
\[
\mathcal A_0(\ell)\xi_0^0(\ell) = -\ell^2\xi_0^1(\ell),\quad 
\mathcal A_0(\ell)\xi_0^1(\ell) = \ell^2\xi_0^0(\ell),
\]
the
$2\times2$ matrix representing the action of $\mathcal A_0(\ell)$ 
on this basis is given by
\[
M_0(\ell) = \begin{pmatrix}0&\ell^2\\-\ell^2&0
\end{pmatrix}.
\]

Next, recall that $\mathcal A_a(\ell)$ anti-commutes with the
reflection $\mathcal S$ defined by \eqref{e:s}. Since 
\[
\mathcal S\xi_0^0(\ell) = -\xi_0^0(\ell),\quad
\mathcal S\xi_0^1(\ell) = \xi_0^1(\ell),
\]
the basis $\{\xi_0^0(\ell),\xi_0^1(\ell)\}$ can be extended to a
basis $\{\xi_a^0(\ell),\xi_a^1(\ell)\}$, for $a\not=0$, with the
same property, i.e., such that $\xi_a^0(\ell)$ and $\xi_a^1(\ell)$
are odd and even in $z$, respectively. In this basis the 
$2\times2$ matrix representing the action of $\mathcal A_a(\ell)$ 
is of the form
\[
M_a(\ell) = \begin{pmatrix}0&\ell^2+O(a)\\-\ell^2+O(a)&0
\end{pmatrix}.
\]

In order to compute the terms of order $a$ we take $\ell=0$.
The derivative $\partial_z P_a$ of the periodic wave always belong to
the kernel of $\mathcal A_a(0)$, 
\[
\mathcal A_a(0)(\partial_z P_a) =0,
\]
due to the invariance of the KP equation under translations in $x$, 
and $\partial_zP_a$ is an odd function, with expansion
\[
\partial_zP_a(z) = -a\sin(z)-\frac16\sin(2z)a^2+O(|a|^3).
\]
Then we choose
\[
\xi_a^0(0) = -\frac1a \partial_zP_a(z) = \sin(z)+\frac16\sin(2z)a+O(a^2),
\]
as a first vector of the basis, which is compatible with the choice
for $a=0$, and satisfies $\mathcal A_a(0)\xi_a^0(0) = 0$. 
For the second vector we set
\[
\xi_a^1(0) =  \cos(z)+ \varphi(z) a+O(a^2),
\]
and 
\[
\mathcal A_a(0)\xi_a^1(0) = (\alpha_0+\alpha_1a+\alpha_2a^2+O(|a|^3))
\xi_a^0(0) .
\]
Then a direct calculation gives
\[
\alpha_0=\alpha_1=0,\quad \alpha_2= -\frac1{12},
\quad \xi_a^1(0) =  \cos(z)+ \frac16\cos(2z) a+O(a^2),
\]
so that 
\[
M_a(0) = \begin{pmatrix}0&-\frac1{12}a^2+O(|a|^3)\\0&0
\end{pmatrix}.
\]
Together with the expression of $M_0(\ell)$, this shows that
\[
M_a(\ell) = \begin{pmatrix}0&\ell^2-\frac1{12}a^2+O(|a|(\ell^2+a^2))
\\-\ell^2 +O(|a|\ell^2)&0
\end{pmatrix}.
\]

The two eigenvalues of $M_a(\ell)$, which are also the
eigenvalues in $\sigma_0(\mathcal A_a(\ell))$, are roots of the
characteristic polynomial
\[
P(\lambda)=\lambda^2+\ell^2(\ell^2-\frac1{12}a^2) + O(|a|\ell^2(\ell^2+a^2)).
\]
Furthermore, since $P_a(z+\pi)=P_{-a}(z)$, the two roots of this
polynomial are the same for $a$ and $-a$, and we conclude that 
\[
\lambda^2 = -\ell^2(\ell^2-\frac1{12}a^2) + O(a^2\ell^2(\ell^2+a^2)).
\]
Consequently, for any $a$ sufficiently small there exists a value
\[
\ell_a^2 = \frac1{12}a^2+O(a^4),
\]
such that
the two eigenvalues are purely imaginary when
$\ell^2>\ell_a^2$, and real, with opposite signs when 
$\ell^2<\ell_a^2$. This proves the lemma.
\end{Proof}

\section{KP-I equation: non-periodic perturbations}
\label{s:5}

In this section we study the transverse spectral stability of the
periodic waves $P_a$ of the KP-I equation, with respect to
perturbations which are localized
or bounded in $z$, i.e.,  we
study the invertibility of the operator 
\[
\mathcal M_a(\lambda,\ell) = \lambda\partial_z-\partial_{z}^4-
\frac1{k_a^2}\partial_{z}^2((1+P_a)\cdot)-\ell^2,
\]
acting in $L^2(\R)$ or $C_b(\R)$,
for $\lambda\in\C$, $\Re\lambda>0$, and $\ell\in\R$, $\ell\not=0$.

\subsection{Bloch-wave decomposition}\label{ss:bloch}

In contrast to the case of perturbations which are periodic in $z$
discussed in the previous section, for perturbations which are
localized or bounded, the linear operator $\mathcal M_a(\lambda,\ell)$ 
acting in $L^2(\R)$ or $C_b(\R)$ has continuous spectrum. A convenient
way to treat this situation is with the help of the so-called
Bloch-wave decomposition, which leads to the following result.

\begin{Lemma}\label{l:bloch}
The linear  operator $\mathcal M_a(\lambda,\ell)$ is invertible, in
either $L^2(\R)$ or $C_b(\R)$, if and only if the linear operators
\[
\mathcal M_a(\lambda,\ell,\gamma) =
\lambda(\partial_z+i\gamma)-(\partial_z+i\gamma)^4- 
\frac1{k_a^2}(\partial_z+i\gamma)^2((1+P_a)\cdot)-\ell^2,
\]
acting in $L^2(0,2\pi)$ with domain $H^4_\per(0,2\pi)$ are invertible,
for any $\gamma\in\left(-\frac12,\frac12\right]$.
\end{Lemma}

\begin{Proof}
The proof relies upon Floquet theory. We refer, e.g., to \cite[Proposition
  1.1 (ii)-(iii)]{bbm} for a detailed proof in similar situation.
\end{Proof}

The key property of the
operators $\mathcal M_a(\lambda,\ell,\gamma)$ in this lemma is that they have
only point spectrum, since they act in $L^2(0,2\pi)$ with
compactly embedded domain $H^4_\per(0,2\pi)$. 
When $\gamma=0$ we recover the case of periodic
perturbations studied in the previous section, so that we restrict now to
$\gamma\not=0$.
Then the operator $\partial_z+i\gamma$ is invertible in $L^2(0,2\pi)$,
property which leads to the following result.

\begin{Lemma}
Assume that $\gamma\in\left(-\frac12,\frac12\right]$ and
  $\gamma\not=0$. Then the linear operator $\mathcal
  M_a(\lambda,\ell,\gamma)$ is invertible in  $L^2(0,2\pi)$ if and
  only if $\lambda$ belongs to the spectrum of the operator
\[
\mathcal A_a(\ell,\gamma) =(\partial_z+i\gamma)^3 +
\frac1{k_a^2}(\partial_z+i\gamma)((1+P_a)\cdot)
+\ell^2(\partial_z+i\gamma)^{-1},
\]
acting in $L^2(0,2\pi)$ with domain $H^3_\per(0,2\pi)$.
\end{Lemma}

We point out that the operator $(\partial_z+i\gamma)^{-1}$ becomes singular,
as $\gamma\to0$, so that the results found by replacing the study of
the invertibility of $\mathcal M_a(\lambda,\ell,\gamma)$ by the study
of the spectrum of $\mathcal A_a(\ell,\gamma)$ are not uniform for
small $\gamma$, i.e., the set of $a$ for which these results are
valid depends upon $\gamma$. We restrict here to 
this latter study which turns out to be simpler, and sufficient to detect
the values of $\ell$ and $\gamma$ where instabilities occur. 

As in the previous section, 
the linear operator $\mathcal A_a(\ell,\gamma) $ decomposes as
\begin{equation}\label{e:dec2}
\mathcal A_a(\ell,\gamma) = 
-(\partial_z+i\gamma)\mathcal L_a(\ell,\gamma),\quad 
\mathcal L_a(\ell,\gamma)=-(\partial_z+i\gamma)^2-
\frac1{k_a^2}\left((1+P_a)\cdot\right)-\ell^2(\partial_z+i\gamma)^{-2}, 
\end{equation}
in which $\partial_z+i\gamma$ and $\mathcal L_a(\ell,\gamma)$ are skew-adjoint
and self-adjoint operators, respectively. As a consequence, the result
in Corollary~\ref{c:imspec} also holds for these operators. 
However, in contrast to the previous
section, here we can also directly apply the result in 
\cite[Theorem 2.13]{HK} showing that the numbers $n(\mathcal
L_a(\ell,\gamma))$ of negative eigenvalues of $\mathcal
L_a(\ell,\gamma)$ (counted with multiplicities)
and $k_u(\mathcal A_a(\ell,\gamma))$ of unstable eigenvalues of 
$\mathcal A_a(\ell,\gamma)$ 
(counted with multiplicities), 
satisfy the inequality
\begin{equation}\label{e:krein}
k_u(\mathcal A_a(\ell,\gamma))\leq n(\mathcal
L_a(\ell,\gamma)),
\end{equation}
provided $\mathcal L_a(\ell,\gamma)$ is invertible.
In fact the result in \cite[Theorem 2.13]{HK} is more precise, but
this inequality is enough for our purposes, since we shall only use this result
in the case $n(\mathcal L_a(\ell,\gamma))=0$. Actually, in this
situation  the
inequality \eqref{e:krein}
also follows from Corollary~\ref{c:imspec}. Indeed, if $n(\mathcal
L_a(\ell,\gamma))=0$ and $\mathcal L_a(\ell,\gamma)$ is invertible,
then the spectrum of $\mathcal L_a(\ell,\gamma)$ is strictly
positive. Since $\mathcal L_a(\ell,\gamma)$ is self-adjoint this
implies the inequality \eqref{e:posL}, so that the spectrum of
$\mathcal A_a(\ell,\gamma)$ is purely imaginary. Consequently,  
$k_u(\mathcal A_a(\ell,\gamma))=0=n(\mathcal
L_a(\ell,\gamma))$,
which proves the inequality.

Another consequence of the decomposition \eqref{e:dec2} is
that the spectrum of $\mathcal
A_a(\ell,\gamma)$ is symmetric with respect to the imaginary
axis, just as the spectrum of $\mathcal A_a(\ell)$ in the previous
section, but it is no longer symmetric with respect to the real axis, 
since the operator $\mathcal A_a(\ell,\gamma)$ is not real. Instead, we
have the following result.

\begin{Lemma}\label{l:spsym2}
Assume that $\gamma\in\left(-\frac12,\frac12\right]$ and
  $\gamma\not=0$. Then the spectrum $\sigma(\mathcal
  A_a(\ell,\gamma))$ of $\mathcal 
A_a(\ell,\gamma)$ is symmetric with respect to the imaginary
axis, and $\sigma(\mathcal A_a(\ell,\gamma))=\sigma(-\mathcal
  A_a(\ell,-\gamma))$. 
\end{Lemma}

\begin{Proof} The first property is a consequence of the decomposition
  \eqref{e:dec2} and of the result in \cite[Proposition 2.5]{HK}. The
  second property is due to the equality
\[
\mathcal A_a(\ell,\gamma) \mathcal S= 
-\mathcal S\mathcal A_a(\ell,-\gamma),
\]
where $\mathcal S$ is the reflection symmetry defined by \eqref{e:s}.
\end{Proof}

\subsection{Transverse spectral instability}

We determine now the spectrum of the operator $\mathcal
A_a(\ell,\gamma)$. 
As a consequence of the second property in Lemma~\ref{l:spsym2} we
restrict 
to $\gamma\in\left(0,\frac12\right]$.
The main result is the following theorem showing the transverse
spectral instability of the periodic waves.

\begin{Theorem}\label{t:nonper}
Assume that $\gamma\in\left(0,\frac12\right]$ and set 
$\ell_c(\gamma) = \sqrt 3\gamma(1-\gamma)$. For any $a$
  sufficiently small, there exists  
$\varepsilon_a(\gamma) = \gamma^{3/2}(1-\gamma)^{3/2}|a|(1+
O(a^2))>0$  such that 
\begin{enumerate}
\item for $|\ell^2-\ell_c^2(\gamma)|\geq
  \varepsilon_a(\gamma)$, the spectrum of $\mathcal
  A_a(\ell,\gamma)$ is purely imaginary;
\item for $|\ell^2-\ell_c^2(\gamma)|<
  \varepsilon_a(\gamma)$, the spectrum of $\mathcal
  A_a(\ell,\gamma)$ is purely imaginary, except for a pair of complex
  eigenvalues with opposite nonzero real parts.
\end{enumerate}
\end{Theorem}

The remainder of this section is occupied by the proof of this theorem.
We rely upon the decomposition \eqref{e:dec2}, the
result in Corollary~\ref{c:imspec} and the inequality \eqref{e:krein}, and
perturbation arguments for linear operators, just as in
Section~\ref{s:4}.

\paragraph{Spectrum of $\boldsymbol{\mathcal L_0(\ell,\gamma)}$ 
and consequences}

\begin{figure}[ht]
\vspace*{3ex}
\begin{center}
\includegraphics[width=3.5cm]{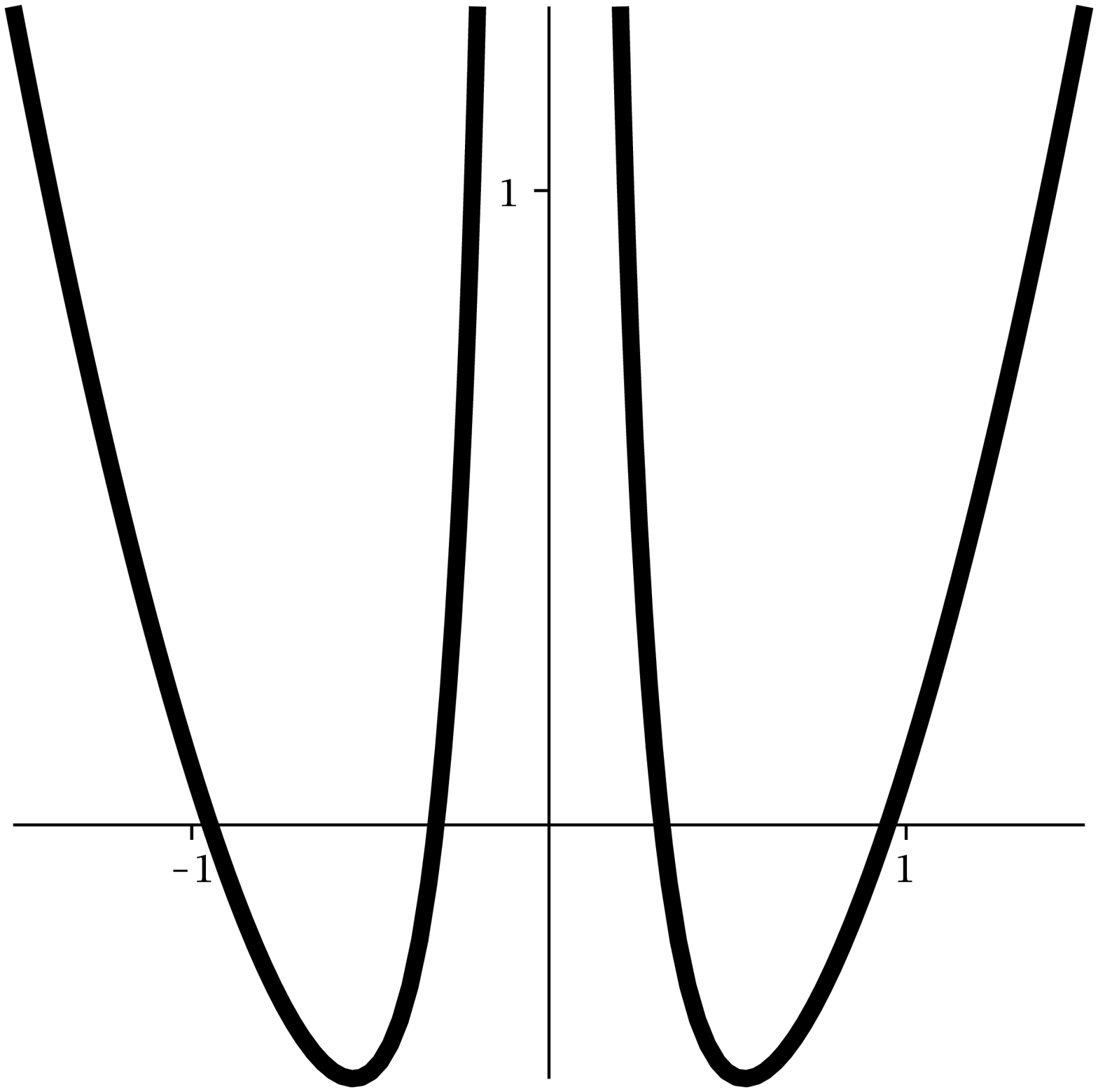}\hspace*{5ex}
\includegraphics[width=3.5cm]{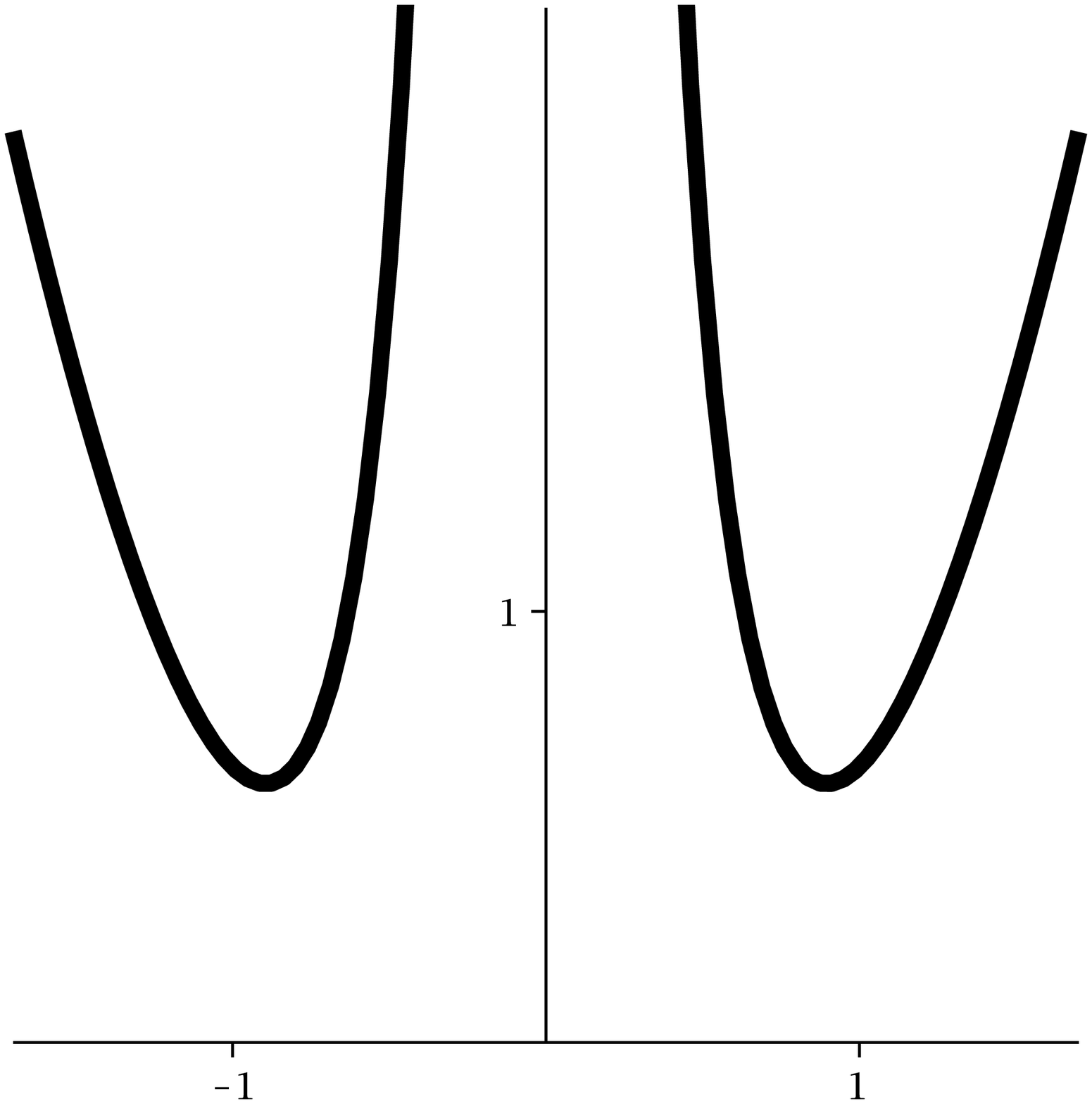}\hspace*{5ex}
\includegraphics[width=3.5cm]{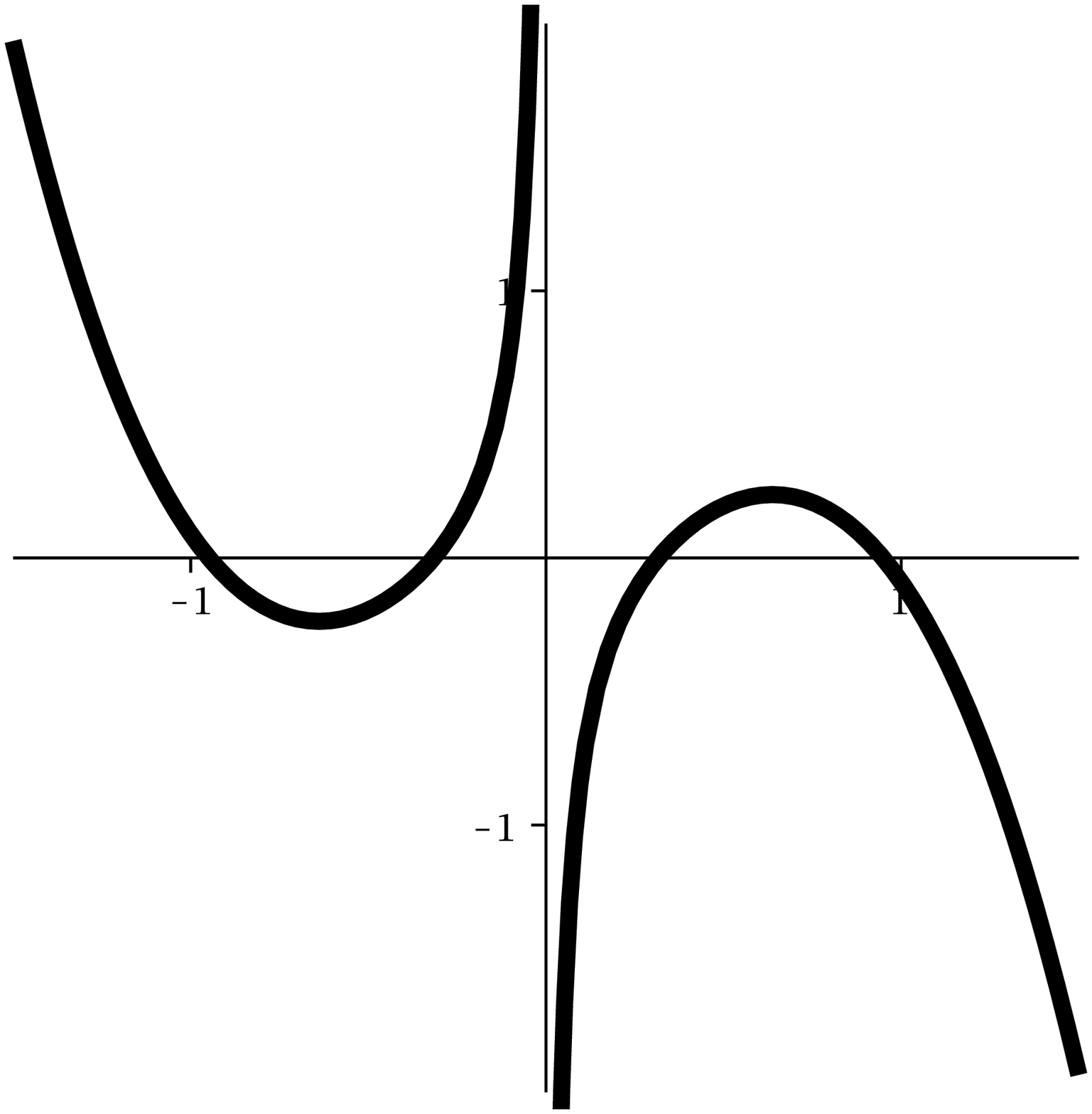}\hspace*{5ex}
\includegraphics[width=3.5cm]{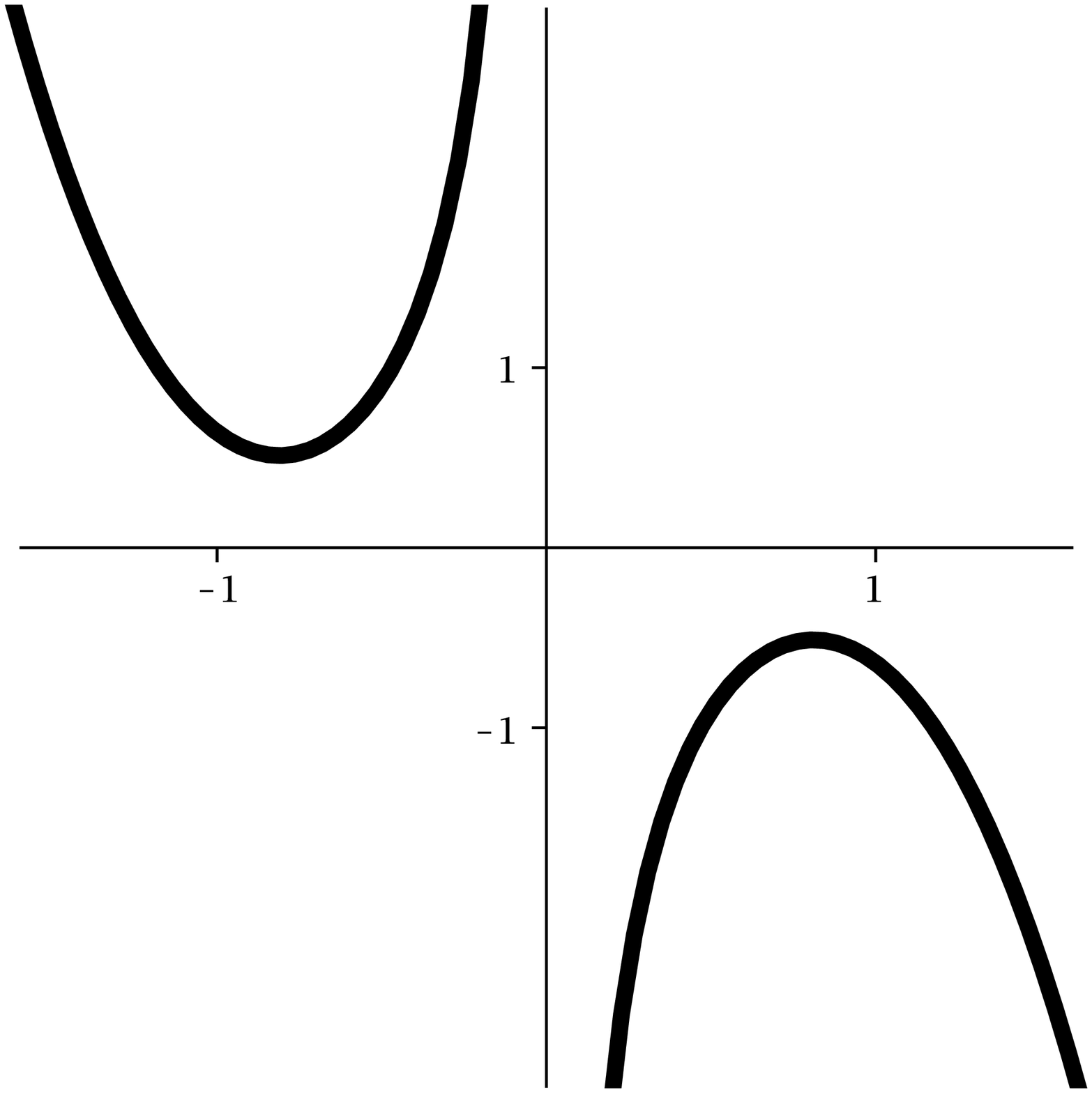}
\end{center}
\centerline{\small\it (a) \hspace*{8.2cm} (b)}
\mcaption{[KP-I equation] \it (a) Graphs of the map $k\mapsto
  k^2-1+\frac{\ell^2}{k^2}$ for  $\ell=0.3$ and $\ell=0.8$ (from left
  to right). The
  eigenvalues of $\mathcal L_0(\ell,\gamma)$ are found by taking
  $k=n+\gamma$, $n\in\Z$. (b) Graphs of the dispersion relation $k\mapsto
  -k^3+k-\frac{\ell^2}{k}$ for the same values of $\ell$. The
  imaginary parts of the 
  eigenvalues of $\mathcal A_0(\ell,\gamma)$ are found by taking
  $k=n+\gamma$, $n\in\Z$. Notice that the zeros of the 
  two maps are the same.}\label{f:disp1}
\end{figure}

The linear operator  $\mathcal L_a(\ell,\gamma)$ is a small
bounded perturbation of 
\[
\mathcal L_0(\ell,\gamma)=-(\partial_z+i\gamma)^2-
1-\ell^2(\partial_z+i\gamma)^{-2}.
\]
Using Fourier series we find  
\[
\sigma(\mathcal L_0(\ell,\gamma))= \left\{
\mu_n(\ell,\gamma)=(n+\gamma)^2-1+\frac{\ell^2}{(n+\gamma)^2}\;;\; n\in\Z
\right\},
\]
which shows that
\begin{itemize}
\item the eigenvalues corresponding to the Fourier modes
  $n\in\Z\setminus\{-1,0\}$ are all  positive,
\[
\mu_n(\ell,\gamma)>(n+\gamma)^2-1\geq \gamma(2+\gamma),\quad\forall\ 
n\in\Z\setminus\{-1,0\};
\]
\item the eigenvalue 
\[
\mu_{-1}(\ell,\gamma)=(1-\gamma)^2-1+\frac{\ell^2}{(1-\gamma)^2}=
\gamma^2-2\gamma+\frac{\ell^2}{(1-\gamma)^2},
\]
corresponding to the Fourier mode $n=-1$ is
  positive when $\ell^2>\ell_-^2$, where $\ell_-^2 =
  \gamma(1-\gamma)^2(2-\gamma)$, it is zero when  $\ell^2=\ell_-^2$,
  and it is negative when  $\ell^2<\ell_-^2$;
\item the eigenvalue
\[
\mu_{0}(\ell,\gamma)=\gamma^2-1+\frac{\ell^2}{\gamma^2},
\] 
corresponding to the Fourier mode $n=0$ is
  positive when $\ell^2>\ell_0^2$, where $\ell_0^2 =
  \gamma^2(1-\gamma^2)$, it is zero when  $\ell^2=\ell_0^2$,
  and it is negative when  $\ell^2<\ell_0^2$;
\end{itemize}
(see also Figure~\ref{f:disp1}~(a)).
Here 
\[
\ell_0^2 =
  \gamma^2(1-\gamma^2) <\ell_-^2 =
  \gamma(1-\gamma)^2(2-\gamma),
\]
for any $\gamma\in\left(0,\frac12\right)$, so that  the
unperturbed operator has positive spectrum for $\ell^2>\ell_-^2$, 
one negative eigenvalue if
  $\ell_0^2<\ell^2<\ell_-^2$, and  two negative eigenvalues if
  $0<\ell^2<\ell_0^2$.
The following result is an immediate consequence of these properties.

\begin{Lemma}
Assume that  $\gamma\in\left(0,\frac12\right]$. For any
  $\varepsilon_*>0$ there exists $a_*>0$, such that the spectrum of  $\mathcal
A_a(\ell,\gamma)$ is purely imaginary,  for any $\ell$ and $a$ satisfying
$\ell^2\geq\ell_-^2+\varepsilon_*$ and 
$|a|\leq a_*$.
\end{Lemma}

\begin{Proof}
The properties above show that 
the eigenvalues of 
the operator $\mathcal L_0(\ell,\gamma)$ are strictly positive,  
$\mu_n(\ell,\gamma)\geq c_0(\gamma)$, for some $c_0(\gamma)>0$,
when $\ell^2\geq\ell_-^2+\varepsilon_*$.
Since $\mathcal L_a(\ell,\gamma)$ is a small
bounded perturbation of  $\mathcal L_0(\ell,\gamma)$  a standard
perturbation argument shows that the eigenvalues of  $\mathcal
L_a(\ell,\gamma)$ remain strictly positive, provided $a$ is
sufficiently small. In particular, $\mathcal L_a(\ell,\gamma)$ is
invertible and $n(\mathcal L_a(\ell,\gamma)=0$, 
so that $k_u(\mathcal A_a(\ell,\gamma))=0$, according to
\eqref{e:krein}, which proves the lemma.
\end{Proof}

It remains to determine the spectrum of $\mathcal A_a(\ell,\gamma)$
for
$\ell^2\in(0,\ell_-^2+\varepsilon_*)$. We proceed as in
Section~\ref{ss:longper}, and decompose the spectrum of  $\mathcal
A_a(\ell,\gamma)$ into $\sigma_0(\mathcal A_a(\ell,\gamma))$
containing a minimal number of eigenvalues, and 
$\sigma_1(\mathcal A_a(\ell,\gamma))$ for which we argue as in
Lemma~\ref{l:long1} to show that it is purely
imaginary. More precisely, $\sigma_1(\mathcal A_a(\ell,\gamma))$ is
such that the restriction of $\mathcal L_a(\ell,\gamma)$ to the
corresponding spectral subspace is positive definite, i.e., it
satisfies the inequality \eqref{e:posL}.

\paragraph{Spectrum of $\boldsymbol{\mathcal A_0(\ell,\gamma)}$ 
and decomposition}

We begin by decomposing the spectrum of the unperturbed operator
$\mathcal A_0(\ell,\gamma)$ in such a way that this decomposition
persists for sufficiently small $a$. 
Using Fourier series again we find  
\[
\sigma(\mathcal A_0(\ell,\gamma))= \left\{
i\omega_n(\ell,\gamma)=
-i(n+\gamma)^3+i(n+\gamma)-i\frac{\ell^2}{n+\gamma}\;;\; n\in\Z
\right\}.
\]
The location on the imaginary axis of these eigenvalues can be viewed
with the help of the dispersion relation
\[
\omega = -k^3+k-\frac{\ell^2}{k};
\]
(see Figure~\ref{f:disp1}~(b)).
Notice that the eigenvalues of $\mathcal A_0(\ell,\gamma)$ which
correspond to negative eigenvalues of $\mathcal L_0(\ell,\gamma)$ are
the ones for $n=-1$, if $\ell_0^2<\ell^2<\ell_-^2$, and for $n=-1$ and
$n=0$, if $0<\ell^2<\ell_0^2$. As a consequence of
Corollary~\ref{c:imspec} there are these eigenvalues which may
lead to instabilities when perturbing the operator  $\mathcal
A_0(\ell,\gamma)$. Therefore, the idea of the spectral decomposition 
is to separate these
eigenvalues from the remaining (infinitely many) eigenvalues, which 
correspond to positive eigenvalues of
$\mathcal L_0(\ell,\gamma)$. 
Of course this separation is possible as long as there are no
collisions with other eigenvalues. Looking for such collisions we find
that
\begin{itemize}
\item for  $\ell_0^2<\ell^2<\ell_-^2$, the eigenvalue
  $i\omega_{-1}(\ell,\gamma)$ collides with
  $i\omega_{0}(\ell,\gamma)$, when
\[
\ell^2=\ell_c^2 = 3\gamma^2(1-\gamma)^2,\quad
\ell_0^2<\ell_c^2<\ell_-^2,\quad\forall\ 
\mbox{$\gamma\in\left(0,\frac12\right)$},
\] 
and it is simple for the other values of $\ell$;
\item for $0<\ell^2\leq\ell_0^2$, the two 
eigenvalues   $i\omega_{-1}(\ell,\gamma)$ and
$i\omega_{0}(\ell,\gamma)$ are simple.
\end{itemize}
This leads to the following result.

\begin{Lemma}\label{l:spA0}
Assume that  $\gamma\in\left(0,\frac12\right]$. There exist
  $\varepsilon_*>0$ and 
  $c_*>0$ such that
\begin{enumerate}
\item  for any $\ell$ satisfying 
$\ell_c^2+\varepsilon_*<\ell^2<\ell_-^2+\varepsilon_*$,
  the spectrum of $\mathcal A_0(\ell,\gamma)$ decomposes as
\[
\sigma(\mathcal A_0(\ell,\gamma))=\{i\omega_{-1}(\ell,\gamma)\}
\cup\sigma_1(\mathcal A_0(\ell,\gamma)),
\]
with
$\mathrm{dist}\big( i\omega_{-1}(\ell,\gamma),
\sigma_1(\mathcal A_0(\ell,\gamma))\big)\geq c_*>0$;
\item  for any $\ell$ satisfying 
$0<\ell\leq\ell_c^2+\varepsilon_*$,
  the spectrum of $\mathcal A_0(\ell,\gamma)$ decomposes as
\[
\sigma(\mathcal A_0(\ell,\gamma))=\{i\omega_{-1}(\ell,\gamma),
i\omega_0(\ell,\gamma)\}
\cup\sigma_1(\mathcal A_0(\ell,\gamma)),
\]
with
$\mathrm{dist}\big(\{i\omega_{-1}(\ell,\gamma),
i\omega_0(\ell,\gamma)\}\,,\,
\sigma_1(\mathcal A_0(\ell,\gamma))\big)\geq c_*>0$.
\end{enumerate}
\end{Lemma}

Perturbation arguments show that this decomposition persists for the
operator $\mathcal A_a(\ell,\gamma)$, for sufficiently small $a$, 
and we argue as in Section~\ref{ss:longper} to locate the spectrum of 
$\mathcal A_a(\ell,\gamma)$.

\paragraph{Spectrum of $\boldsymbol{\mathcal A_a(\ell,\gamma)}$ for 
 $\boldsymbol{\ell_c^2+\varepsilon_*<\ell^2<\ell_-^2+\varepsilon_*}$}
We prove the following result showing that the spectrum of  $\mathcal
  A_a(\ell,\gamma)$ is purely imaginary for such values of $\ell$, 
provided  $a$ is sufficiently small.

\begin{Lemma}\label{l:lc1}
Assume that  $\gamma\in\left(0,\frac12\right]$. There exist
$\varepsilon_*>0$ and $a_*>0$ such that the spectrum of $\mathcal
  A_a(\ell,\gamma)$ is purely imaginary, for any $\ell$ and $a$
  satisfying $\ell_c^2+\varepsilon_*<\ell^2<\ell_-^2+\varepsilon_*$
  and $|a|\leq a_*$.
\end{Lemma}

\begin{Proof}
The operator $\mathcal  A_a(\ell,\gamma)$ is a small relatively
bounded perturbation of the operator  $\mathcal  A_0(\ell,\gamma)$,
and there exists  $c_1>0$ such that
\[
\| \mathcal  A_a(\ell,\gamma) - \mathcal  A_0(\ell,\gamma)\|_{H^1\to
  L^2}\leq c_1|a|, 
\]
for any  $\ell>0$ and
$a$ sufficiently small. Then by arguing as in the proof of
Lemma~\ref{l:spdec}, taking into account 
the spectral decomposition in Lemma~\ref{l:spA0}, we conclude that the
spectrum of  $\mathcal  A_a(\ell,\gamma)$ decomposes as
\[
\sigma(\mathcal  A_a(\ell,\gamma)) = \sigma_0(\mathcal
A_a(\ell,\gamma)) \cup \sigma_1(\mathcal  A_a(\ell,\gamma)),
\]
where $\sigma_0(\mathcal  A_a(\ell,\gamma))$ and $\sigma_1(\mathcal
A_a(\ell,\gamma))$ have the following properties
\begin{enumerate}
\item there exist two open balls $B_0\subset B_1$ centered at
  $i\omega_{-1}(\ell,\gamma)$ such that
\[
\sigma_0(\mathcal A_a(\ell,\gamma))\subset B_0,\quad
\sigma_1(\mathcal A_a(\ell,\gamma))\subset \C\setminus\overline{B_1};
\]
\item the spectral projection $\Pi_a(\ell,\gamma)$ associated with 
$\sigma_0(\mathcal A_a(\ell,\gamma)) $ satisfies
  $\|\Pi_a(\ell,\gamma)-\Pi_0(\ell,\gamma)\|=O(|a|)$;
\item the spectral subspace $\mathcal X_a(\ell,\gamma) =
  \Pi_a(\ell)(L^2_0(0,2\pi))$ is one dimensional;
\end{enumerate}
(see also  Lemma~\ref{l:spdec}). This decomposition holds when
  $\ell_c^2+\varepsilon_*<\ell^2<\ell_-^2+\varepsilon_*$ 
  and $|a|\leq a_*$, for $\varepsilon_*$ and $a_*$ are sufficiently
  small. 

Next, the eigenvalues in $\sigma_1(\mathcal A_0(\ell,\gamma))$
correspond to positive eigenvalues of $\mathcal L_0(\ell,\gamma))$, so that 
we can argue as in the proof of Lemma~\ref{l:long1} to conclude that
the eigenvalues in $\sigma_1(\mathcal A_a(\ell,\gamma))$ are purely
imaginary. Finally, the property (iii) shows that  
$\sigma_0(\mathcal A_a(\ell,\gamma))$ consists of one simple eigenvalue. This
eigenvalue is the unique eigenvalue of $\mathcal A_a(\ell,\gamma)$ in
the open ball $B_0$, and it is purely imaginary 
since the spectrum of  $\mathcal
A_a(\ell,\gamma)$ is symmetric with respect to the imaginary
axis. This completes the proof of the lemma.
\end{Proof}

\paragraph{Spectrum of $\boldsymbol{\mathcal A_a(\ell,\gamma)}$ for 
 $\boldsymbol{0<\ell^2\leq\ell_c^2+\varepsilon_*}$}
Finally, we prove the following result showing that there is a pair of
unstable eigenvalues of $\mathcal A_a(\ell,\gamma)$, for $\ell$ close
to the value $\ell_c$ at which 
the two eigenvalues of $\mathcal A_0(\ell,\gamma)$
corresponding to the Fourier modes $n=-1$ and $n=0$ collide.
This result completes the proof of Theorem~\ref{t:nonper}.

\begin{Lemma}\label{l:lc2}
Assume that  $\gamma\in\left(0,\frac12\right]$. There exist
positive constants  $\varepsilon_*$,  $a_*$, and 
$\varepsilon_a(\gamma) = \gamma^{3/2}(1-\gamma)^{3/2}|a|(1+
O(a^2))\leq\varepsilon_*$  such that the spectrum of $\mathcal
  A_a(\ell,\gamma)$ is purely imaginary, if
$0<\ell^2\leq\ell_c^2+\varepsilon_*$ and $|\ell^2-\ell_c^2|\geq
  \varepsilon_a(\gamma)$. If  $|\ell^2-\ell_c^2|<
  \varepsilon_a(\gamma)$ then  the spectrum of $\mathcal
  A_a(\ell,\gamma)$ is purely imaginary, except for a pair of complex
  eigenvalues with opposite nonzero real parts.
\end{Lemma}

\begin{Proof}
As in the previous case, using 
the spectral decomposition for $a=0$ in Lemma~\ref{l:spA0}, and the
arguments in the proof of  Lemmas~\ref{l:spdec} and \ref{l:long1}
we conclude that the
spectrum of  $\mathcal  A_a(\ell,\gamma)$ decomposes as
\[
\sigma(\mathcal  A_a(\ell,\gamma)) = \sigma_0(\mathcal
A_a(\ell,\gamma)) \cup \sigma_1(\mathcal  A_a(\ell,\gamma)),
\]
where $\sigma_1(\mathcal  A_a(\ell,\gamma))$ is purely imaginary, and 
$\sigma_0(\mathcal  A_a(\ell,\gamma))$ consists now of two eigenvalues
which are the continuation for small $a$ of the eigenvalues
$i\omega_{-1}(\ell,\gamma)$ and $i\omega_0(\ell,\gamma)$.
This decomposition holds for any $\ell$ and $a$ satisfying
$0<\ell^2\leq\ell_c^2+\varepsilon_*$ and $|a|\leq a_*$, with 
$\varepsilon_*$ and $a_*$ sufficiently small.
Consequently, it remains to locate the two eigenvalues in 
$\sigma_0(\mathcal  A_a(\ell,\gamma))$.

The two eigenvalues $i\omega_{-1}(\ell,\gamma)$ and
$i\omega_0(\ell,\gamma)$ are simple and the distance between them is
strictly positive
\[
|i\omega_{-1}(\ell,\gamma)-i\omega_0(\ell,\gamma)|\geq c_0,
\]
for any $\ell$ outside a neighborhood of $\ell_c$, where the
eigenvalues collide, and such that
$0<\ell^2\leq\ell_c^2+\varepsilon_*$.
A standard perturbation argument then shows that the continuation of
this eigenvalues for sufficiently small $a$ is a pair of simple
eigenvalues of $\mathcal  A_a(\ell,\gamma)$. Each of these eigenvalues
of $\mathcal  A_a(\ell,\gamma)$ is purely imaginary, since the
spectrum of $\mathcal  A_a(\ell,\gamma)$ is symmetric with respect to the
imaginary axis, just as in the proof of Lemma~\ref{l:lc1}. This shows
that the spectrum of $\mathcal  A_a(\ell,\gamma)$ is purely imaginary
for any $\ell$ outside a neighborhood of $\ell_c$.

In order to locate $\sigma_0(\mathcal  A_a(\ell,\gamma))$ for $\ell$
close to $\ell_c$, we 
proceed as in the proof of
Lemma~\ref{l:long0}. We compute successively
a basis for the two-dimensional spectral subspace associated with 
$\sigma_0(\mathcal  A_a(\ell,\gamma))$, 
the $2\times2$ matrix $M_a(\ell,\gamma)$ 
representing the action of $\mathcal A_a(\ell,\gamma)$
on this basis, and the eigenvalues of this matrix.

At $a=0$ we choose as basis the two eigenvectors associated with the
eigenvalues $i\omega_0(\ell,\gamma)$ and  $i\omega_{-1}(\ell,\gamma)$,
\[
\xi_0^0(\ell,\gamma) = 1,\quad \xi_0^1(\ell,\gamma)  =e^{-iz},
\]
and then the
$2\times2$ matrix representing the action of $\mathcal A_0(\ell,\gamma)$ 
on this basis is given by
\[
M_0(\ell,\gamma) = \begin{pmatrix}i\omega_0(\ell,\gamma)&0
\\0&i\omega_{-1}(\ell,\gamma)
\end{pmatrix}.
\]
Next, notice that 
the operator $\mathcal A_a(\ell,\gamma)$ anti-commutes with the
conjugation-reflection symmetry $\mathcal S_c$ defined by 
\[
\mathcal S_c(u(z))=\overline{u(-z)}.
\]
Since 
\[
\mathcal S_c\xi_0^0(\ell,\gamma) = \xi_0^0(\ell,\gamma),\quad
\mathcal S_c\xi_0^1(\ell,\gamma) = \xi_0^1(\ell,\gamma),
\]
the basis $\{\xi_0^0(\ell,\gamma),\xi_0^1(\ell,\gamma)\}$ can be extended to a
basis $\{\xi_a^0(\ell,\gamma),\xi_a^1(\ell,\gamma)\}$, for $a\not=0$, with the
same property. In this basis, the
$2\times2$ matrix $M_a(\ell,\gamma)$ anti-commutes with the
conjugation (the map induced by $\mathcal S_c$), i.e.,
$M_a(\ell,\gamma)=-\overline{M_a(\ell,\gamma)}$, 
which implies that the coefficients of  
$M_a(\ell,\gamma)$ are purely imaginary. 

In order to compute the terms of order $a$, we take $\ell=\ell_c$, and
proceed as in the
computation of the vector $\xi_a^1(0)$ in the proof of
Lemma~\ref{l:long0}. Then at order $a$ we find
\[
M_a(\ell_c,\gamma) =
\begin{pmatrix}i\omega_0(\ell_c,\gamma)&
\frac i2\gamma a
\\
\frac i2(\gamma-1) a&i\omega_{-1}(\ell_c,\gamma)
\end{pmatrix}+O(a^2).
\]
Together with the expression of $M_0(\ell,\gamma)$, this shows that
\[
M_a(\ell,\gamma) = 
\begin{pmatrix}i\omega_0(\ell_c,\gamma)-i\frac\varepsilon\gamma&
\frac i2\gamma a
\\
\frac i2(\gamma-1) a&i\omega_{-1}(\ell_c,\gamma)-i\frac\varepsilon{\gamma-1}
\end{pmatrix}+O(|a|(|\varepsilon|+|a|)),
\]
where $\varepsilon=\ell^2-\ell_c^2$.

Recall that
\[
\ell_c^2 = 3\gamma^2(1-\gamma)^2,\quad 
i\omega_0(\ell_c,\gamma)=i\omega_{-1}(\ell_c,\gamma)=
-2i\gamma(1-\gamma)(1-2\gamma). 
\]
Then looking for eigenvalues $\lambda$ of the form
\[
\lambda = -2i\gamma(1-\gamma)(1-2\gamma) +iX,
\]
we find that $X$ is root of the polynomial
\[
P(X)=X^2+X\left(\frac\varepsilon\gamma - \frac\varepsilon{1-\gamma} +
O(a^2)\right)
+\frac14\gamma(1-\gamma)a^2-\frac{\varepsilon^2}{\gamma(1-\gamma)} +
  O(a^2(|\varepsilon|+a^2)) .
\]
This polynomial has real coefficients, 
since the matrix $M_a(\ell,\gamma)$ is purely imaginary. Here
we have also taken in account the fact that 
these coefficients are even in $a$, due to the equality
$P_a(z+\pi)=P_{-a}(z)$, which implies that  the two roots of the
polynomial are the same for $a$ and $-a$.
A direct computation shows that the discriminant of this polynomial is
\[
\Delta_a(\varepsilon,\gamma) =
\frac{\varepsilon^2}{\gamma^2(1-\gamma)^2} -\gamma(1-\gamma)a^2+
O(a^2(|\varepsilon|+a^2)), 
\]
and for any $a$ sufficiently small there exists
\[
\varepsilon_a(\gamma) = \gamma^{3/2}(1-\gamma)^{3/2}|a|(1+
O(a^2))>0,
\]
such that
the two roots of the polynomial are real when
$|\varepsilon|\geq\varepsilon_a(\gamma)$, 
and complex otherwise. This implies that the two eigenvalues of
$\mathcal A_a(\ell,\gamma)$ are purely imaginary when 
$|\ell^2-\ell_c^2|\geq\varepsilon_a(\gamma)$, and complex, with
opposite nonzero real parts when $|\ell^2-\ell_c^2|<\varepsilon_a(\gamma)$,
which proves the lemma.
\end{Proof}

\section{KP-II equation: periodic perturbations}
\label{s:6}

In this section we briefly discuss the  
KP-II equation, when $\sigma=-1$ in equation \eqref{e:kp}.
The functional set-up developed in
Section~\ref{ss:specper} for periodic perturbations, and in
Section~\ref{ss:bloch} for non-periodic perturbations remains valid
for the KP-II equation. However, it turns out that locating the
spectra of the resulting operators is more complicated, 
because of the qualitatively different behavior of the
dispersion relation. We shall restrict here to the case
of periodic perturbations, when we have to
determine the spectrum of the operator
\[
\mathcal A_a(\ell) = -\partial_z\mathcal L_a(\ell),\quad 
\mathcal L_a(\ell)=-\partial_{z}^2-
\frac1{k_a^2}\left((1+P_a)\cdot\right)+\ell^2\partial_z^{-2}, 
\]
acting in $L^2_0(0,2\pi)$ with
domain $H^3_\per(0,2\pi)\cap L^2_0(0,2\pi)$ (see Corollary~\ref{c:inv}
in Section~\ref{ss:specper}).
We begin by analyzing the spectra of the unperturbed operators
$\mathcal L_0(\ell)$ and $\mathcal A_0(\ell)$.

\begin{figure}[ht]
\hspace*{3ex}
\begin{center}
\includegraphics[width=3.5cm]{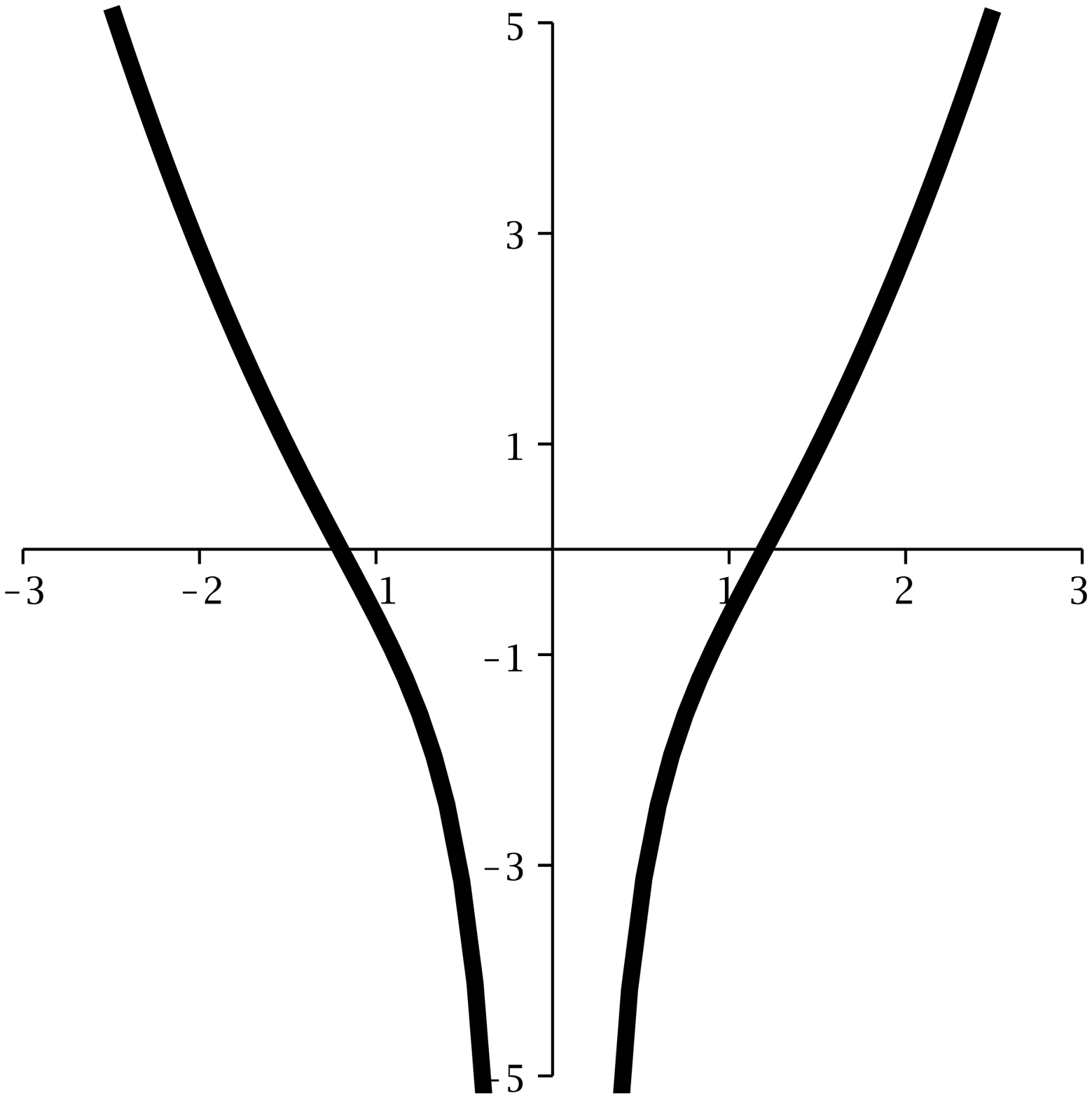}\hspace*{10ex}
\includegraphics[width=3.5cm]{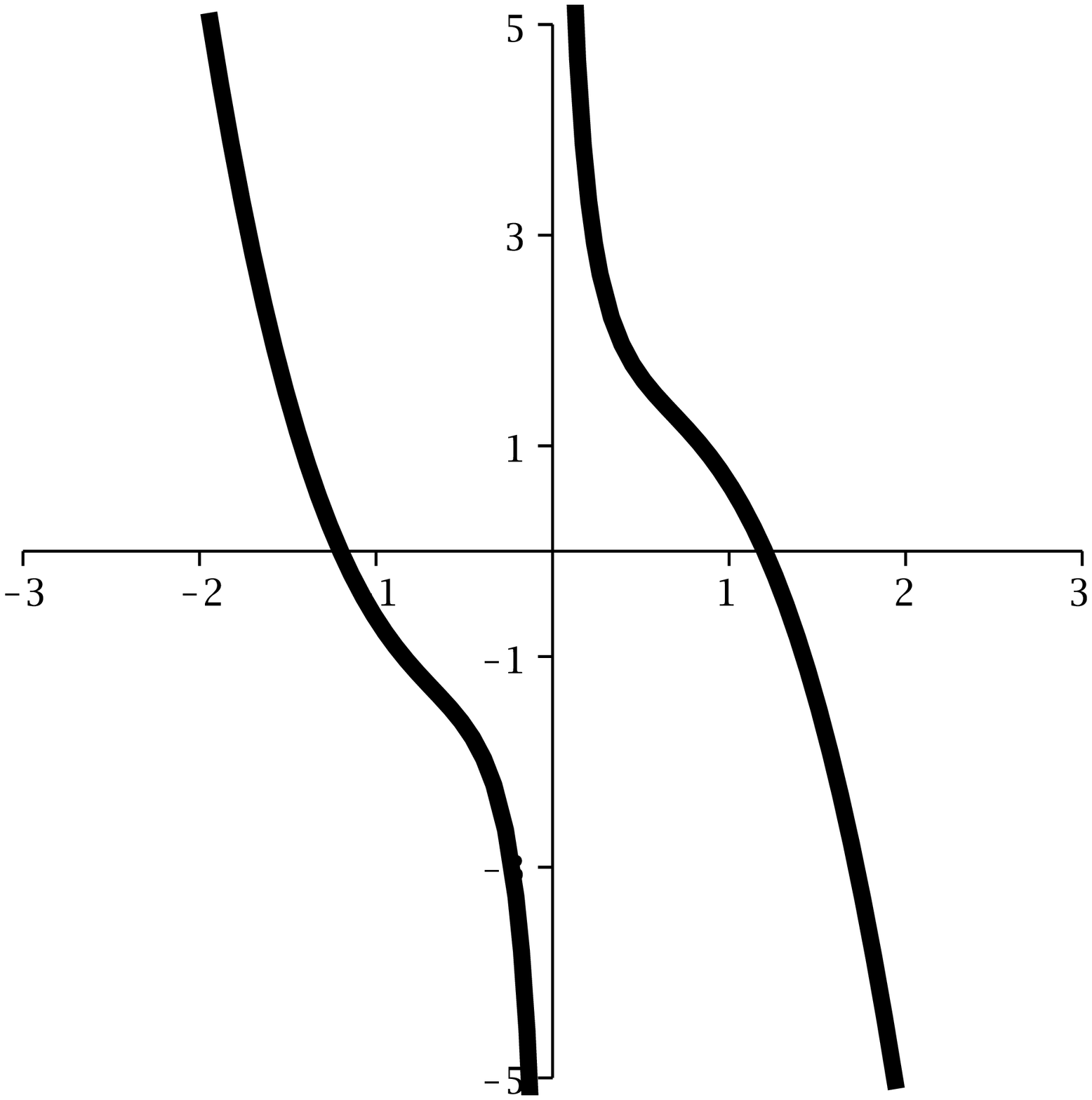}
\end{center}
\centerline{\small\it (a) \hspace*{4.6cm} (b)}
\mcaption{[KP-II equation] \it (a) Graph of the map $k\mapsto
  k^2-1-\frac{\ell^2}{k^2}$ for $\ell=0.8$. The
  eigenvalues of $\mathcal L_0(\ell)$ are found by taking
  $k=n$, $n\in\Z^*$. (b) Graph of the dispersion relation $k\mapsto
  -k^3+k+\frac{\ell^2}{k}$ for  $\ell=0.8$. The
  imaginary parts of the 
  eigenvalues of $\mathcal A_0(\ell)$ are found by taking
  $k=n$, $n\in\Z^*$. Notice that the zeros of the 
  two maps are the same.}\label{f:disp2}
\end{figure}

\paragraph{Spectrum of $\boldsymbol{\mathcal L_0(\ell)}$}
Using Fourier series we find that the spectrum of the operator 
$\mathcal L_0(\ell)$
acting in  $L^2_0(0,2\pi)$ is given by
\[
\sigma(\mathcal L_0(\ell))=
\left\{\mu_n(\ell)=n^2-1-\frac{\ell^2}{n^2}\;;\;n\in\Z^*\right\};
\]
(see also Figure~\ref{f:disp2} (a)). 
In contrast to the KP-I equation,
the spectrum of $\mathcal L_0(\ell)$ contains now negative
eigenvalues, and the number of these eigenvalues increases with
$\ell$. Therefore we cannot use the result in Corollary~\ref{c:imspec}
to exclude unstable eigenvalues of $\mathcal A_a(\ell)$ for certain
values of $\ell$, e.g., as in Lemma~\ref{l:short}. However, for a
given $\ell$ we can use this corollary to reduce the number of
potentially unstable eigenvalues to a finite number,  e.g., as in
Lemmas~\ref{l:long1}, \ref{l:lc1}, or \ref{l:lc2}, but this number tends
to $\infty$ as $\ell\to\infty$.

\paragraph{Spectrum of $\boldsymbol{\mathcal A_0(\ell)}$}
The spectrum of the operator 
$\mathcal A_0(\ell)$
acting in  $L^2_0(0,2\pi)$ is given by
\[
\sigma(\mathcal A_0(\ell))=
\left\{i\omega_n(\ell)=-in^3+in+i\frac{\ell^2}{n}\;;\;n\in\Z^*\right\};
\]
(see also Figure~\ref{f:disp2} (b)). Notice that the dispersion
relation 
\begin{equation}\label{e:disp2}
\omega = -k^3+k+\frac{\ell^2}k
\end{equation}
is monotonically decreasing on $(-\infty,-1]$ and $[1,\infty)$, so
that colliding eigenvalues correspond to Fourier modes with opposite
signs. A direct calculation then shows that for any $m,p\in\N^*$ the
eigenvalues corresponding to the Fourier modes $m$ and $-p$ collide
when
\[
\ell^2=\ell_{m,p}^2 = mp(m^2-mp+p^2-1).
\]
Moreover, the corresponding eigenvalues of $\mathcal L_0(\ell)$ have
opposite signs, so that any of
these collisions may lead to unstable eigenvalues of the operator 
$\mathcal A_a(\ell)$. 
We can use the arguments presented in the previous sections to reduce
the spectral analysis to that of the continuation, for small $a$, of
the two colliding eigenvalues, but an additional difficulty comes now
from the fact that there are infinitely many values
of $\ell$ where two potentially unstable eigenvalues collide.  
Furthermore, it turns out that in order to
determine the location of these two eigenvalues,
we need to compute
Taylor expansions up to order $m+p$ for
the resulting $2\times2$ matrix. 
Consequently, the complete spectral analysis seems to be much more
difficult in this case. However, we can easily conclude in the case
$m=p=1$, which corresponds to values of  $\ell$ close to $\ell_{1,1}=0$.

\paragraph{Long wavelength transverse perturbations}

Restricting to small
values of $\ell$, i.e., long wavelength transverse
perturbations,
we can argue as in Section~\ref{ss:longper} and determine
the spectrum of $\mathcal A_a(\ell)$. 
In contrast to the KP-I equation, now 
we show that the spectrum of  $\mathcal A_a(\ell)$ is
purely imaginary when $\ell$
and $a$ sufficiently small, i.e., the small periodic waves of the KP-II
equation are spectrally stable with respect to long
wavelength transverse perturbations. 

\begin{Lemma}
The spectrum  of  $\mathcal A_a(\ell)$ is purely
imaginary, for any $\ell$ and $a$ sufficiently small. 
\end{Lemma}

\begin{Proof}
It is enough to replace $\ell^2$ by $-\ell^2$ in the proofs given  
in Section~\ref{ss:longper}. Then the results in
Lemmas~\ref{l:spdec} and \ref{l:long1} remain the same, the only
change appears in the result of the calculation of the eigenvalues of 
the $2\times2$ matrix $M_a(\ell)$ in
Lemma~\ref{l:long0}. Upon  replacing $\ell^2$ by $-\ell^2$
we find that the two eigenvalues satisfy
\[
\lambda^2 = -\ell^2(\ell^2+\frac1{12}a^2) + O(a^2\ell^2(\ell^2+a^2)).
\]
Consequently, $\lambda^2<0$, for $\ell$ and $a$ sufficiently
small, which shows that the two eigenvalues are purely imaginary, 
so that the spectrum of  $\mathcal A_a(\ell)$ is purely
imaginary.
\end{Proof}

\section{Discussion}
\label{s:7}

We have presented a functional set-up for the study of the transverse
spectral stability of the one-dimensional periodic waves of the KP
equation. Relying upon perturbation arguments for linear operators we
have determined the spectra of the operators arising in the spectral
problems for the waves
of small amplitude of the KP-I equation, and showed that these waves
are transversely unstable with respect to perturbations which are either
periodic or non-periodic (localized or bounded) in the direction of
propagation. For the KP-II equation, the study turns out
to be more complicated, due to the `negative dispersion', which leads
to a great number of `potentially unstable' modes. In this case
it is only shown that the periodic waves are transversely spectrally stable 
with respect to perturbations which are periodic in the direction of
propagation and have long wavelengths in the transverse
direction. It seems much more difficult to study 
transverse perturbations in the finite and short wavelength
regimes, with this type of arguments, 
since the location of the potentially unstable eigenvalues is
determined by terms of arbitrary high order in the expansion of these
eigenvalues for small $a$. In these regimes 
a numerical computation of the spectrum, using for instance the 
package described in \cite{DCK,DK}, may help to locate the spectrum,
and so get more insight in the stability properties of these waves.

This type of approach to transverse spectral stability of periodic waves can 
be adapted to other two-dimensional dispersive models, e.g., 
generalized KP equations, Boussinesq
systems, equations of Schr\"odinger type. 
However, one has to restrict to periodic
waves of small amplitude, since we strongly rely upon perturbation
arguments for linear operators. For large waves, we mention the approach in
\cite{JZ} which leads to instability criteria, and in the particular
case of integrable
equations, so also for the KP equation,  
one could take advantage of the integrability and explicitly compute
the spectra of the linear operators (e.g., see \cite{BD} for such an
approach in the case of the KdV
equation). 

Finally, we point out that these results
concern only the spectral
stability of the periodic waves, the question of nonlinear 
stability being widely open. Nevertheless, we expect that the recent
methods developed in  \cite{RT1,RT2} for the transverse nonlinear
instability of the solitary waves, can be
adapted to periodic waves.
In particular, these would lead to a transverse nonlinear instability
result for the periodic waves of the KP-I equation.

\end{document}